\documentclass[11pt]{article} 

\renewcommand{\arraystretch}{1.2}
\title{Identifiability and estimation of recursive max-linear models}

\author{Nadine Gissibl\thanks{Center for Mathematical Sciences, Technical University of Munich,  85748 Garching, Boltzmannstrasse 3, Germany; e-mail: n.gissibl@tum.de, cklu@tum.de}
	\and
	Claudia Kl\"uppelberg\footnotemark[1]
	\and
	Steffen Lauritzen\thanks{Department of Mathematical Sciences, University of Copenhagen, Universitetsparken 5, 2100 Copenhagen, Denmark; e-mail: lauritzen@math.ku.dk}
}

\usepackage{lmodern}
\usepackage{siunitx}
\usepackage{color}
\usepackage[usenames, dvipsnames]{xcolor}
\usepackage{verbatim}
\usepackage{amssymb}
\usepackage{amsfonts}
\usepackage{amsmath}
\usepackage{mathabx}
\usepackage{amsthm}
\usepackage[numbers]{natbib}
\usepackage[small,bf,figurewithin=section]{caption}
\usepackage{graphicx}
\usepackage{graphics}
\usepackage{scrpage2}
\usepackage{geometry}
\usepackage{nicefrac}
\usepackage{multicol}
\usepackage{multirow}
\usepackage{enumerate}
\usepackage{relsize}
\usepackage{dsfont}
\usepackage{bbm}
\usepackage{url}
\usepackage{capt-of}
\usepackage[english]{babel}
\usepackage[]{inputenc}
\usepackage{accents}
\usepackage{pgfplots}
\usepackage{tikz}
\usetikzlibrary{shapes.misc,positioning,arrows,shapes}
\tikzset{cross/.style={cross out, draw=black, minimum size=2*(#1-\pgflinewidth), inner sep=0pt, outer sep=0pt},
	cross/.default={1pt}}
\usepackage{chngcntr}
\usepackage{float}
\usepackage{subfig}
\usepackage[symbol]{footmisc}
\usepackage{paralist}
\usepackage{tocvsec2}
\usepackage{pdflscape} 
\usepackage{longtable}
\usepackage{capt-of}
\usepackage{pgfplots}
\usepackage{afterpage}
\usepackage{rotating}

\numberwithin{equation}{section}

\theoremstyle{plain}
\newtheorem{theorem}{Theorem}[section]
\newtheorem{lemma}[theorem]{Lemma}
\newtheorem{remark}[theorem]{Remark}
\newtheorem{proposition}[theorem]{Proposition}
\newtheorem{definition}[theorem]{Definition}
\newtheorem{corollary}[theorem]{Corollary}
\newtheorem{example}[theorem]{Example}

\newtheorem{algorithm}[theorem]{Algorithm}

\newcommand{\bthe}{\begin{theorem}}
	\newcommand{\ethe}{\end{theorem}}

\newcommand{\ble}{\begin{lemma}}
	\newcommand{\ele}{\end{lemma}}

\newcommand{\bde}{\begin{definition}\rm}
	\newcommand{\ede}{\end{definition}}

\newcommand{\bco}{\begin{corollary}}
	\newcommand{\eco}{\end{corollary}}

\newcommand{\bpr}{\begin{proposition}}
	\newcommand{\epr}{\end{proposition}}

\newcommand{\brem}{\begin{remark}\rm}
	\newcommand{\erem}{\end{remark}}

\newcommand{\bproof}{\begin{proof}}
	\newcommand{\eproof}{\end{proof}}

\newcommand{\bexam}{\begin{example}\rm}
	\newcommand{\eexam}{\end{example}}

\newcommand{\balg}{\begin{algorithm}\rm}
	\newcommand{\ealg}{\end{algorithm}}

\def\calb{{\mathcal{B}}}
\def\calm{{\mathcal{M}}}
\def\call{{\mathcal{L}}}

\def\R{{\mathbb R}}

\def\P{{\mathbb P}}
\def\borel{{\mathbb {B} }}

\newcommand{\bfz}{\boldsymbol{Z}}

\newcommand{\bfx}{\boldsymbol{X}}
\newcommand{\bfsx}{\boldsymbol{x}}

\newcommand{\bfy}{\boldsymbol{Y}}
\newcommand{\bfsy}{\boldsymbol{y}}

\def\bone{\mathds{1}}
\def\1{\mathds{1}}

\def\D{\mathcal{D}}
\def\calp{{\mathcal{P}}}
\def\calf{{\mathcal{F}}}

\newcommand{\DAG}{\text{DAG}}

\newcommand{\mSEM}{\text{recursive ML model}}

\newcommand{\ML}{\text{ML}}

\newcommand{\CM}{{\text{ML coefficient matrix}}}

\newcommand{\ubpp}{{\rm{supp}}}
\newcommand{\an}{{\rm{an}}}
\newcommand{\pa}{{\rm{pa}}}

\newcommand{\des}{{\rm{de}}}
\newcommand{\An}{{\rm{An}}}
\newcommand{\Pa}{{\rm{Pa}}}

\newcommand{\wh}{\widehat}
\newcommand{\wt}{\widetilde}

\newcommand{\halmos}{\quad\hfill\mbox{$\Box$}}  

\definecolor{plum}{cmyk}{0.50,1,0,0}
\definecolor{TealBlue}{cmyk}{0.86,0,0.34,0.02}
\definecolor{OliveGreen}{cmyk}{0.64,0,0.95,0.40}

\newcommand{\Nadine}{\textcolor{red}}

\allowdisplaybreaks[4]

\begin{document}
\maketitle

\subsection*{Abstract}
We address the identifiablity and estimation of recursive max-linear structural equation models represented by an edge weighted directed acyclic graph (\DAG). Such models are generally unidentifiable and 
we identify the whole class of \DAG s and edge weights corresponding to a given observational distribution. 
For estimation, standard likelihood theory cannot be applied because the corresponding families of distributions are not dominated.  Given the  underlying \DAG, we present  an estimator for the class of edge weights and show that it can be considered a  generalized maximum likelihood estimator. In addition, we develop a simple method for identifying the structure of the DAG. With probability tending to one at an exponential rate with the number of observations, this method correctly identifies the class of DAGs and, similarly, exactly identifies the possible edge weights.
\bigbreak
\noindent {\em MSC 2010 subject classifications:} Primary 60E15, 62H12; secondary 62G05, 60G70, 62-09
\bigbreak
\noindent {\em Keywords and phrases:} Causal inference, Bayesian network, directed acyclic graph, extreme value theory, generalized maximum
likelihood estimation,  graphical model, identifiability, max-linear model, 
structural equation model.

\section{Introduction} 

Establishing and understanding cause-effect relations is an omnipresent desire in science and daily life. It is especially important when dealing with  extreme events, because they are mostly dangerous and very costly; knowing and understanding the causes of such events and their causal relations could help us to deal better with them. Examples include incidents at airplane landings (\citet{LH1}), 
flooding in river networks (\citet{engelke2015,engelke:hitz:18}), financial risk (\citet{einmahl2016}), and  chemical pollution of rivers (\citet{Hoef2006}).  Such applications, where extreme risks may propagate through a network, have been the motivation behind the definition of  {\em recursive max-linear\/} (ML) models in \citet{GK1}. Recursive \ML\ models are {\em structural equation models\/} (SEMs)  represented by a \emph{directed acyclic graph} (\DAG) and thereby obey the basic Markov properties associated with directed graphical models (\citet{lauritzen:01,Lauritzen1990}).  Both SEMs (see for example \citet{Bollen}, \citet{Pearl2009}) and directed graphical models (see for example \citet{KF,Lauritzen1996,SGS}) are well-established concepts for the understanding and quantification of causal inference from observational data.  We note that \citet{hitz:evans:16} and \citet{engelke:hitz:18} discuss graphical models for extremes that are based on undirected graphs.

Recursive \ML\ models are defined by a \DAG, a collection of edge weights, and a vector of independent innovations.  Important research problems that are addressed for recursive SEMs are the question of {\em identifiability} of the coefficients and the associated \DAG\ from the observational distribution.  
Although the true \DAG\ and edge weights for a recursive \ML\ model are not identifiable from the observational distribution, the so-called {\em max-linear coefficient matrix\/} is identifiable and determines the possible class of \DAG s and edge weights uniquely.

We shall show that estimation and structure learning of recursive \ML\ models   can be done in a simple and efficient fashion 
by exploiting properties of the ratios between observable components of the model.
For a sufficiently large number of observations, these ratios identify the true  \ML\ coefficient matrix with a probability that converges exponentially fast to 1. 
For the situation where the DAG is known, we show that our estimator can be considered a maximum likelihood estimator in an extended sense, originally introduced by \citet{kiefer1956}.

\medskip

Our paper is organized as follows. In Section~\ref{ch4:s2} we introduce the model class of recursive \ML\ models and  the notation used throughout. In Section~\ref{ch4:s4} we discuss the identifiability of a recursive \ML\ model from its observational distribution. Here we show distributional properties of the ratio between two components. 
Based on these properties, we suggest an identification method. Section~\ref{ch4:s5} is then devoted to the estimation of recursive \ML\ models where we assume  the \DAG\ to be known. We show that the proposed estimates are  {\em generalized maximum likelihood estimates  (GMLEs)} in the sense of  Kiefer--Wolfowitz.  The main part is here the derivation of a specific Radon-Nikodym derivative. 
In  Section~\ref{ch4:s6} we complement the theoretical findings on the identifiability of recursive ML models with an efficient procedure to learn recursive \ML\ models from observations only, even when the DAG itself is also unknown. 
Section~\ref{ch4:s7} concludes and suggests further directions of research.

\section{Preliminaries --- recursive max-linear models\label{ch4:s2}} 

In this section we introduce notation and summarize the most important properties of recursive \ML\ models needed.
A \mSEM\ for  a random vector $\bfx=(X_1,\ldots,X_d)$ is specified by an underlying structure in terms of a \DAG\ $\D$ with nodes $V=\{1,\ldots, d\}$, positive {\em edge weights\/} $c_{ki}$ for $i\in V$ and $k\in \pa(i)$, and independent positive random variables $Z_1, \dots, Z_d$  with support $\R_+:=(0,\infty)$ 
and atom-free distributions:
\begin{align}\label{ch4:ml-sem}
X_i=\bigvee_{k \in {{\pa}}(i)} c_{ki} X_k \vee  Z_i, \quad i=1,\ldots,d,
\end{align}
where $\pa(i)$ are the parents of node $i$ in $\D$. To highlight the \DAG\ $\D$, we say that $\bfx$ follows a  \emph{\mSEM} on $\D$. Note that this is a slight variation of the original definition in \cite{GK1}. We shall refer to 
 $\bfz=(Z_1,\ldots, Z_d)$ 
 as the vector of \emph{innovations}.

In the context of risk analysis, natural candidates for distributions of the innovations are extreme value distributions or distributions in their domain of attraction, resulting in a corresponding multivariate distribution (for details and background on multivariate extreme value models, see for example \citet{Beirlant2004,   DHF,Resnick1987,Resnick2007}).  

Throughout the paper we use the following notation. The sets $\an(i)$, $\pa(i)$, and $\des(i)$ 
contain the ancestors, parents, and  descendants  of node $i$ in $\D$. 
We set $\An(i)=\an(i)\cup\{i\}$ and $\Pa(i)=\pa(i)\cup\{i\}$. For $U\subsetneq V$  we write $\bfx_U=(X_\ell, \ell\in U)$ and accordingly for $\bfsx\in \R_+^d$, $\bfsx_U=(x_\ell, \ell\in U)$. 

Instead of $k\in\pa(i)$ we also write $k\to i$. Assigning the weight $d_{ji}(p)= \prod_{\nu=0}^{n-1} c_{k_\nu k_{\nu+1}}$ to every path $p=[j=k_0 \to k_1 \to \dots \to k_n=i]$ and denoting the set of all paths from $j$ to $i$ by $P_{ji}$, the  non-negative matrix $B=(b_{ij})_{d\times d}$ with entries  
\begin{align}\label{ch4:coeff}
b_{ji}=\bigvee_{p\in P_{ji}}d_{ji}(p) \quad\text{for } j\in \an(i), \quad b_{ii}=1,\quad  \text{and}  \quad b_{ji}=0 \quad\text{for }  j\in V\setminus \An(i),
\end{align}
is said to be the {\em \ML\ coefficient matrix of $\bfx$}. This means for distinct $i,j\in V$, $b_{ji}$ is positive if and only if there is a path from $j$ to $i$; in that case $b_{ji}$ is the maximum weight of all paths from $j$ to $i$, where the weight of a path is the product of all edge weights $c_{ki}$ along this path. We say that a path from $j$ to $i$ whose weight equals $b_{ji}$ is {\em max-weighted}. 

The components of $\bfx$ can also be expressed as max-linear functions of their ancestral innovations and an independent one; the corresponding {\em \ML\ coefficients} are the entries of $B$:
\begin{align}\label{ch4:ml-noise}
X_i=\bigvee_{j=1}^db_{ji}Z_j=\bigvee_{j\in\An(i)}b_{ji}Z_j, \quad i=1,\ldots, d;
\end{align} 
see Theorem~2.2 of \cite{GK1}.

For two non-negative matrices $F$ and $G$, where the number of columns in $F$ is equal to the number of rows in $G$, we define the matrix product $\odot: \overline\R_+^{m\times n}\times \overline\R_+^{n\times p} \rightarrow \overline\R_+^{m\times p}$
by
\begin{align}\label{ch2:odot} 
(F=(f_{ij})_{m\times n},G=(g_{ij})_{n\times p}) \mapsto F\odot G :=\Big(\bigvee\limits_{k=1}^n f_{ik}g_{kj}\Big)_{m\times p},
\end{align}where $\overline\R_+=[0,\infty)$.
The triple  $(\overline \R_+,\vee,\cdot)$, 
is an idempotent semiring with $0$ as 0-element and  $1$ as 1-element and 
the operation $\odot$ is therefore a matrix product over this semiring; see for example\ \citet{But2010}. 
Denoting by $\calm$ all $d\times d$ matrices with non-negative entries and by $\vee$ the componentwise maximum between two matrices, $(\calm,\vee,\odot)$ is also a semiring with the null matrix as 0-element and the $d\times d$ identity matrix  $I_d$ 
as 1-element.  

The matrix product $\odot$ allows us to represent the ML coefficient matrix $B$ of $\bfx$
 in terms of the weighted adjacency matrix $(c_{ij} \bone_{\pa(j)}(i))_{d\times d}$ of $\D$ since \eqref{ch4:coeff} and \eqref{ch4:ml-noise} simply become
 \begin{align}\label{eq:dotrepr}
 B= (I_d\vee C)^{\odot (d-1)}= \bigvee_{k=0}^{d-1}  C^{\odot k}, \qquad \bfx= \bfz\odot B,
  \end{align}
where we have let $A^{\odot 0}=I_d$ and $A^{\odot k}= A^{\odot(k-1)}\odot A$ for $A\in \overline\R_+^{d\times d}$ and $k\in\mathbb{N}$; see Proposition 1.6.15 of \citet{But2010} as well as Theorem~2.4 and Corollary 2.5 of \cite{GK1}. 
\medskip

\section{Identifiability of a recursive max-linear model\label{ch4:s4}} 

In this section we discuss the question of identifiability of the elements of a recursive \ML\ model  from the distribution $\call(\bfx)$ of $\bfx$.  Indeed we shall show the following:

\bthe\label{ch4:identify}
Let $\call(\bfx)$ be the distribution of $\bfx$ following a recursive \ML\ model. 
 Then its \ML\ coefficient matrix $B$ and the distribution of its innovation vector $\bfz$ are identifiable from $\call(\bfx)$. Furthermore, the class of all \DAG s and edge weights that could have generated $\bfx$ by \eqref{ch4:ml-sem} can be obtained.
\ethe  

The remaining part of this section is devoted to proving Theorem~\ref{ch4:identify}, but first we shall consider a small  example, illustrating the issues. 

\bexam[The DAG and the edge weights are not necessarily identifiable\label{ch4:exprob2}]\\
Consider a \mSEM\  
on the \DAG\ $\D$ depicted below with edge weights $c_{12}, c_{23}, c_{13}$.
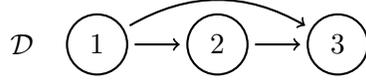
\begin{center}
		\begin{tikzpicture}[->,every node/.style={circle,draw},line width=0.8pt, node distance=1.6cm,minimum size=0.8cm,outer sep=1mm]
		\node (1) [outer sep=1mm]  {$1$};
		\node (2) [right of=1,outer sep=1mm] {$2$};
		\node (3) [right of=2,outer sep=1mm] {$3$}; 
		\foreach \from/\to in {2/3}
		\draw (\from) -- (\to);   
		\foreach \from/\to in {1/2}
		\draw (\from) -- (\to);   
		\foreach \from/\to in {1/3}
		\draw [ bend left]  (1) to (3);
		\node (n5)[draw=white,fill=white,left of=1,node distance=1cm] {$\D$};
		\end{tikzpicture}
\end{center}According to \eqref{ch4:ml-sem}, the components of $\bfx$ have the following representations
\begin{align*}
X_1 = Z_1,\quad X_2 = c_{12}X_1\vee Z_2,\quad \text{and} \quad  X_3 = c_{13}X_1 \vee c_{23}X_2 \vee Z_3.
\end{align*}
but also representations in terms of the innovations using \eqref{ch4:ml-noise} as
\begin{align*}
X_1 = Z_1,\quad X_2 = c_{12}Z_1\vee Z_2,\quad \text{and} \quad  X_3 = ( c_{12}c_{23}\vee c_{13})Z_1 \vee c_{23}Z_2 \vee Z_3,
\end{align*}
If  $c_{13} \le c_{12}c_{23}$ we have for any $c^*_{13} \in [0,c_{12}c_{23}]$
that $b_{13}= c_{12}c_{23}\vee c^*_{13}=c_{12}c_{23}\vee c_{13}=c_{12}c_{23}$; so 
 we could also write
\begin{align*}
X_3=c^*_{13}X_1\vee c_{23}X_2\vee Z_3
\end{align*}
without changing the distribution $\call(\bfx)$ of $\bfx$.
This implies that if $c_{13}\leq c_{12}c_{23}$,  $\bfx$  follows a recursive ML model on $\D$ with edge weights $c_{12},c_{23},c^*_{13}$ but it also follows a recursive model on the \DAG\ $\D^B$ depicted below with edge weights $c_{12}, c_{23}$. 
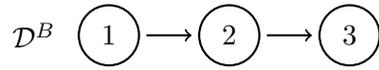
\begin{center}
	\begin{tikzpicture}[->,every node/.style={circle,draw},line width=0.8pt, node distance=1.6cm,minimum size=0.8cm,outer sep=1mm]
	\node (1) [outer sep=1mm]  {$1$};
	\node (2) [right of=1,outer sep=1mm] {$2$};
	\node (3) [right of=2,outer sep=1mm] {$3$}; 
	\foreach \from/\to in {2/3}
	\draw (\from) -- (\to);   
	\foreach \from/\to in {1/2}
	\draw (\from) -- (\to);  
	\node (n5)[draw=white,fill=white,left of=1,node distance=1cm] {$\D^B$}; 
	\end{tikzpicture}
\end{center}

Consequently, we can neither identify $\D$ nor the value $c_{13}$ from the distribution $\call(\bfx)$ of $\bfx$. However, note that the \ML\ coefficient  $b_{13}=c_{12}c_{23}\vee c_{13}$ is uniquely determined. 
If we however assume that $c_{13}>c_{12}c_{23}$, only $\D$ 
and  the edge weights $c_{12}, c_{23}, c_{13}$ represent $\bfx$ in the sense of \eqref{ch4:ml-sem}. Thus in this case the \DAG\ and the edge weights are identifiable from the distribution $\call(\bfx)$. 
\halmos
\eexam

As conclusion of Example~\ref{ch4:exprob2}, it is generally not possible to identify the true  \DAG\  $\D$ and  the edge weights $c_{ki}$  underlying $\bfx$ in representation \eqref{ch4:ml-sem} from $\call(\bfx)$, since several \DAG s and edge weights  may exist such that   $\bfx$ has this  representation. The smallest \DAG\ of this kind 
is the \DAG\ that has an edge $k\to i$ if and only if $k\to i$ is the only max-weighted path from $k$ to $i$. We call this \DAG\ $\D^B$ the {\em minimum \ML\ \DAG\ of $\bfx$} and note that this is uniquely determined from the ML coefficient matrix $B$. 
All other \DAG s representing 
$\bfx$ are those that include the edges of $\D^B$ and whose nodes have the same ancestors. The edge weights $c_{ki}$ in the representation \eqref{ch4:ml-sem} of $\bfx$ are only uniquely determined for edges contained in $\D^B$;  namely, by $b_{ki}$; otherwise, $c_{ki}$ may be any number in $(0,b_{ki}]$.  We summarize these findings in the following theorem which is paraphrasing 
Theorems~5.3 and 5.4 of \cite{GK1}.
\bthe \label{thm:paraphrase} Suppose $\bfx$ follows a recursive ML model with edge weights $C=\{c_{ij}\}$ and ML coefficient matrix $B$. Let $\D^B$ be the minimum ML DAG of $\bfx$ as described above.  Then a DAG  $\D^*$ with associated weight matrix $C^*$ is a valid representation of $\bfx$ if and only if
\begin{enumerate}
\item $\D^B \subseteq \D^*$;
\item $\D^*$ and $\D^B$ have the  same reachability matrix;
\item $c^*_{ij}=c_{ij}$ for $i\in \pa^B(j)$;
\item $c^*_{ij}\in (0,b_{ij}]$ for $i\in \pa^*(j)\setminus \pa^B(j)$,
\end{enumerate}
where  $\pa^B(j)$ and $\pa^*(j)$ denote the parents of $j$ in $\D^B$ and $\D^*$ respectively.  
\ethe

Based on the above observations, we investigate the identifiability of the whole class  of \DAG s and edge weights representing the max-linear structural equations \eqref{ch4:ml-sem} of $\bfx$ from  $\call(\bfx)$.  Since this class can be recovered from $B$, it suffices to clarify whether $B$ is identifiable from   $\call(\bfx)$.   There are many ways to prove that this is indeed the case.  The way we present in this section  suggests a simple procedure to estimate $B$  from independent realizations of $\bfx$ (see Algorithm~\ref{ch4:alg4} below). An alternative way can be found in Appendix~4.A.1 of \cite{Gissibl2018}.

The ratios $\bfy=\{Y_{ij}=X_j/X_i,\; i,j=1\ldots, d\}$ 
between all pairs of components of $\mathbf{X}$ are the essential quantities used to identify $B$ from $\call(\bfx)$. We first present distributional properties of these ratios, 
 where we let  $(\Omega,\mathcal{F},\P)$ denote  the probability space of $\bfz$  
 and, hence,  of $\bfx$. In what follows, we use the standard convention and write  events such as  $\{\omega\in\Omega: X_i(\omega)< X_j(\omega) \}$ 
  as $\{X_i<X_j\}$, etc. 
Unsurprisingly, because of the max-linear representation \eqref{ch4:ml-noise} of the components of $\bfx$, the ratios inherit their distributional properties from the innovations. 
It plays an important role that 
\begin{align}\label{ch4:noiseNull}
\text{the event $\big\{Z_i=x Z_j\big\}$ for distinct  $i,j\in V$ and $x\in \R_+$  has probability zero,}
\end{align}
which follows from  the independence of the innovations and the fact that their distributions are atom-free.

\ble\label{lem:XioverXj} Let $i,j\in V$ be distinct.
\begin{enumerate}
	\item[(a)] The ratio $Y_{ji}=X_i/X_j$ 
	has an atom in $x\in\mathbb{R}_+$ 
	if and only if $\An(i)\cap\An(j)\neq \emptyset$ and $x={b_{\ell i}}/{b_{\ell j}}$ for some $\ell\in\An(i)\cap \An(j)$.
	\item[(b)] We have
	\begin{align*}
	\ubpp(Y_{ji})=
	\begin{cases}
[b_{ji},\infty) 
&\text{if $j\in\an(i)$}\\
\big( 0,{1}/{b_{ij}}\big] 
&\text{if $j\in\des(i)$} \\
\mathbb{R}_+ 	
& \text{otherwise}, 
	\end{cases}
	\end{align*}
	where $\ubpp(Y_{ji})$ denotes the support of $Y_{ji}$. 
\end{enumerate}
\ele
\bproof To establish 
(a) note that \eqref{ch4:ml-noise} and \eqref{ch4:noiseNull} imply that  the sets 
$ \{X_i=xX_j\}=\big\{ \bigvee_{\ell\in \An(i)} b_{\ell i}Z_\ell= \bigvee_{\ell\in \An(j)} x b_{\ell j}Z_\ell\big\} $ and 
\begin{align*}
\Big\{\bigvee_{\ell \in\An(i)\cap \An(j): \atop b_{\ell i}= b_{\ell j}x} b_{\ell i}Z_\ell > \bigvee_{\ell\in \An(i)\cap \An(j):\atop  b_{\ell i}\neq  b_{\ell j}x } (b_{\ell i} \vee x b_{\ell j}) Z_\ell \vee \bigvee_{\ell\in\An(i)\setminus \An(j)} b_{\ell i} Z_\ell \vee  \bigvee_{\ell\in \An(j)\setminus \An(i)} xb_{\ell j} Z_\ell  \Big\} 
\end{align*} 
differ only by a set of probability zero.  Since  the innovations are independent and have support $\R_+$ 
 the conclusion follows. \\

To establish (b) note that the support $\R_+$ of the innovations   
and the representation  \eqref{ch4:ml-noise} yield
\begin{align*}
\ubpp(Y_{ji})= \left\{ \frac{\bigvee_{\ell\in \An(i)} b_{\ell i}z_\ell }{\bigvee_{\ell\in \An(j)} b_{\ell j}z_\ell }: \boldsymbol{z}_{\An(i)\cup\An(j)} \in \R_+^{\vert \An(i)\cup\An(j)\vert} \right\}.
\end{align*}
The continuity of the function 
\begin{align*}
\R_+^{\vert \An(i)\cup \An(j)\vert} \rightarrow \R_+, \quad \boldsymbol{z}_{ \An(i)\cup \An(j)} \mapsto \frac{\bigvee_{\ell\in \An(i)} b_{\ell i}z_\ell }{\bigvee_{\ell\in \An(j)} b_{\ell j}z_\ell }
\end{align*}
implies that $ \ubpp(Y_{ji})$ is an interval in $\R_+$.  Since for $j\in\text{an}(i)$ by Corollary~3.13 of \cite{GK1}  $b_{ji} \le Y_{ji}$  and by (a) $b_{ji}$ is an atom of $Y_{ji}$, it suffices to show that $j\in\text{an}(i)$ if  $ \ubpp(Y_{ji})$ has a positive lower bound. For this assume that $j\not\in\an(i)$. Because of the positive lower bound of $ \ubpp(Y_{ji})$, there exists some  $a\in \R_+$ 
such that
\begin{align}\label{ch4:lowboundsupp}
\bigvee_{\ell\in \An(i)\cap \An(j)}a b_{\ell j}z_\ell\vee \bigvee_{\ell\in \An(j)\setminus \An(i)} ab_{\ell j}z_\ell  \le \bigvee_{\ell\in \An(i) } b_{\ell i}z_\ell 
\end{align}
for all $ \boldsymbol{z}_{\An(i)\cup \An(j)}\in \R_+^{\vert\An(i)\cup\An(j)\vert}$. 
As $\An(j)\setminus \An(i)\neq\emptyset$,  for fixed $\mathbf{z}_{\text{An}(i)}\in\R_+^{\vert \An(i)\vert}$, we can choose $z_\ell$ for some $\ell\in \An(j)\setminus \An(i)$ so large that $ab_{\ell j}z_\ell$ is greater than  the maximum on the right-hand side of \eqref{ch4:lowboundsupp}. This contradicts  \eqref{ch4:lowboundsupp}. Hence, $j\in\an(i)$.
\eproof

In Table~\ref{ch4:table:rel} we summarize the results of Lemma~\ref{lem:XioverXj}: depending on the relationship between $i$ and $j$ in $\D$, the support and atoms of $Y_{ji}$  are shown. 

\begin{table}[htb]
	\begin{center}
		\caption{Distributional properties of $Y_{ji}$ for distinct $i,j \in V$.\label{ch4:table:rel}}
		\begin{tabular}{l|c|l}
			Relationship between $i$ and $j$&   $\ubpp(Y_{ji})$ & Atoms  \\
			\hline
			$j\in\an(i)$&  $[b_{ji},\infty)$ 
			& $\{{b_{\ell i}}/{b_{\ell j}}, \ell\in \An(j)\} $ \\
			$i\in\an(j)$&$(0,{1}/{b_{ij}}]$ 
			& $\{{b_{\ell i}}/{b_{\ell j}}, \ell\in \An(i)\}$ \\
			otherwise:
			& \\
			\hspace*{1em} if $\an(i)\cap \an(j)\neq \emptyset$ & $\R_+$ 
			&  $\{{b_{\ell i}}/{b_{\ell j}}, \ell\in \an(i)\cap \an(j)\}$ \\
			\hspace*{1em} if $\an(i)\cap \an(j)= \emptyset$ &$\R_+$ 
			&$\emptyset$
		\end{tabular}
	\end{center}
\end{table}

Table~\ref{ch4:table:rel} and the fact that $b_{ji}=0$ for $j\not\in\An(i)$ (cf.~\eqref{ch4:coeff}) suggest the following algorithm to find $B$ from $\call(\bfx)$ since  we can identify the support of $Y_{ji}$ from $\call(\bfx)$. This proves  the identifiability of $B$ from $\call(\bfx)$.  In fact, it is sufficient to know  $\ubpp(Y_{ji})$ for all $i,j\in V$ with $i\neq j$ rather than the whole distribution $\call(\bfx)$.

\balg[Find $B$ from $\call(\bfx)$\label{ch4:alg1}] 
\begin{enumerate}
	\item[1.] For all $i\in V=\{1,\ldots,d\}$, set $b_{ii}=1$. 
	\item[2.] For all  $i,j \in V$ with $i\neq  j$,  find $\ubpp(Y_{ji})$:
	\item[] if  $\ubpp(Y_{ji})=[a,\infty)$ for some $a\in \R_+$, then set $b_{ji}=a$;
	\begin{enumerate}
		\item[] else, set $b_{ji}=0$.   
	\end{enumerate}
\end{enumerate}
\ealg

So far we have shown that the  \CM\  $B$ of $\bfx$ can be obtained from $\call(\bfx)$. Since all \DAG s and edge weights that represent $\bfx$ in the sense of \eqref{ch4:ml-sem} can be determined from $B$, the only quantities  we do not know about yet but appear in the definition of $\bfx$ are the innovations. In what follows we show that the distribution of the innovation vector $\bfz$ is also identifiable from $\call(\bfx)$. For this, due to the identifiability of $B$ from $\call(\bfx)$ and the independence of the innovations, it suffices to provide an algorithm that  determines the distributions of the innovations from $\call(\bfx)$ and $B$.  Note that $B$ also determines the ancestral relationships between any pair of nodes in that $j\in \An(i)$ for any DAG representing $\bfx$ if and only if $b_{ji}>0$.

We denote by $F_{Z_i}$ the distribution function of the innovation $Z_i$. For this algorithm, we do not have to know the whole distribution $\call(\bfx)$; it is enough to know the ML coefficient matrix $B$ and the univariate marginal distribution functions of $\call(\bfx)$.

\balg[Find $F_{Z_1}(x),\ldots, F_{Z_d}(x)$ for  $x\in\R_+$ 
from $B$ and $\call(\bfx)$\label{ch4:alg3}]\\
For $\nu=0,\ldots, d-1$,\\
\hspace*{2em} for $i\in V$ such that $\vert\an(i)\vert=\vert \{j\in V\setminus\{i\}: b_{ji}\neq 0\}\vert=\nu$, set
\begin{align*}
F_{Z_i}(x) =\frac{\P(X_i\le x)}{\prod_{j\in\an(i)} F_{Z_j}({x}/{b_{j i}})}. 
\end{align*}
\ealg
Here we have used the convention that $\prod_{j\in\emptyset} a_j=1.$ The correctness of Algorithm~\ref{ch4:alg3} follows from the independence of the innovations and representation \eqref{ch4:ml-noise}.

\section{Estimation with known directed acyclic graph \label{ch4:s5}} 

In this section we consider independent realizations $\bfsx^{(t)}=\big(x^{(t)}_1,\ldots,x^{(t)}_d\big)$,    $t=1,\ldots, n$, of a random vector $\bfx=(X_1,\ldots,X_d)$ following a recursive \ML\ model  with its \DAG\ $\D$ given. 
Further, we consider the distribution of the innovation vector
to be fixed; however, we emphasize that our estimates and their validity do not depend on this distribution as long as it prescribes 
independent, atom-free margins with support $\R_+$. 
Our aim is the  estimation of the edge weights $c_{ki}$ and the \ML\ coefficient matrix $B$.
We recall from Theorem~\ref{thm:paraphrase} that only the \ML\ coefficient matrix $B$ can be directly identified from $\call(\bfx)$ and hence our focus will be on the estimation of $B$;  subsequently all DAGs and systems of edge weights compatible with $B$ can be obtained from Theorem~\ref{thm:paraphrase}.

\subsection*{The ML coefficient matrix $\boldsymbol{B}$} 

In the following we let $\mathcal{B}(\D)$ denote  the class of possible 
\ML\ coefficient matrices of all recursive \ML\ models on $\D$.  
For $B$ being a matrix with non-negative entries and diagonal elements $b_{ii}=1$ we define $B_0:= (b_{ij} \bone_{\pa(j)}(i))_{d\times d}$. Then it holds that  $B\in \mathcal{B}(\D)$ if and only if $B$ satisfies the following
\begin{equation}\label{eq:fixpoint}[b_{ji}>0 \iff j\in \An(i)] \text { and }B=I_d\vee (B\odot B_0); 
\end{equation}
 see Theorem~4.2 or   Corollary~4.3(a) of \cite{GK1}. 
 
\bexam[Illustration of (\ref{eq:fixpoint})]
To illustrate the above, consider the small network below
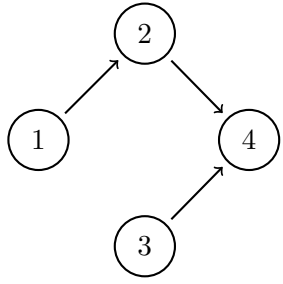
\begin{center}
\begin{tikzpicture}[->,every node/.style={circle,draw},line width=0.8pt, node distance=1 cm,minimum size=0.8cm,outer sep=1mm]

      \node(a) at (0,0) {$1$};
      \node(b) [above right = of a] {$2$};
      \node(c) [below right = of a] {$3$}; 
      \node(d) [right =18mm of a]{$4$};

		\draw (a) -- (b);
   \draw (c) -- (d);
	\draw (b) -- (d);

    \end{tikzpicture}
    \end{center}and  a potential ML coefficient matrix $B$ with reduction $B_0$, as given below. 
    \begin{align*}
B = 
\begin{pmatrix}
b_{11} &  b_{12}    & 0 & b_{14} \\
0 & b_{22} &   0    &  b_{24} \\
0   &    0   &   b_{33}&   b_{34} \\
0 &0  &0 &  b_{44}
\end{pmatrix} \quad B_0 = 
\begin{pmatrix}
0 &  b_{12}    & 0 & 0 \\
0 & 0 &   0    &  b_{24} \\
0   &    0   &   0&   b_{34} \\
0 &0  &0 &  0
\end{pmatrix}, 
\end{align*}
where we have used that $1$ and $2$ are not ancestors of $3$ and $1$ is not a parent of $4$. We wish to check whether $B\in B(\D)$ for this particular DAG so we further calculate
\[ I_4\vee (B\odot B_0)= I_4\vee \begin{pmatrix}
0 &  b_{11}b_{12}    & 0 & b_{12}b_{24} \\
0 & 0 &   0    &  b_{22}b_{24} \\
0   &    0   &   0&   b_{33}b_{34} \\
0 &0  &0 &  0
\end{pmatrix} =  \begin{pmatrix}
1 &  b_{11}b_{12}    & 0 & b_{12}b_{24} \\
0 & 1 &   0    &  b_{22}b_{24} \\
0   &    0  &   1&   b_{33}b_{34} \\
0 &0  &0 &  1
\end{pmatrix}\]
Now $B=I_4\vee (B\odot B_0)$ readily implies that  $b_{ii}=1, i=1,\ldots,4$ and $b_{14}=b_{12}b_{24}$.
\eexam
\subsection*{A simple estimate of $\boldsymbol{B}$\label{ch4:s52}} 

Next we discuss  a sensible estimate  of $B$. 
Table~\ref{ch4:table:rel} shows that for $j\in\an(i)$ the minimal value that can be observed for the ratio $Y_{ji}=X_i/X_j$ is $b_{ji}$, which is an atom of 
$Y_{ji}$. This suggests 
the following estimate $\breve B$ of the \ML\ coefficient matrix: 
\begin{align*}
\breve  b_{ii}=1,\,\,\,  \breve b_{ji}=0 \text{ for }  j\in V\setminus \An(i), \,\,\, \text{and}\,\,\,\breve  b_{ji}=\bigwedge_{t=1}^n y^{(t)}_{ji}  = \bigwedge_{t=1}^n \frac{x_i^{(t)}}{x_j^{(t)}}   \text{ for } j\in \an(i). 
\end{align*}
\citet{Davis1989} suggested such minimal observed ratios as estimates for parameters in max-ARMA processes. 
For $n$  sufficiently large, we can  expect  to observe the atoms $b_{ji}$ for $j\in\an(i)$ in the sample $\bfsx^{(1)}, \ldots , \bfsx^{(n)}$ and, hence, to estimate  the \ML\ coefficients exactly.
 However, if $n$ is not large we may with positive probability have that $\breve B$ is not an 
 \ML\ coefficient matrix of any recursive \ML\ model on $\D$ as the following simple example shows:

 \bexam[$\breve B$ is not necessarily in $\calb(\D)$]\\
Consider the \DAG
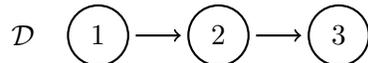
\begin{center}
\begin{tikzpicture}[->,every node/.style={circle,draw},line width=0.8pt, node distance=1.6cm,minimum size=0.8cm,outer sep=1mm]
 \node (1) [outer sep=1mm]  {$1$};
  \node (2) [right of=1,outer sep=1mm] {$2$};
  \node (3) [right of=2,outer sep=1mm] {$3$}; 
        \foreach \from/\to in {2/3}
  \draw (\from) -- (\to);   
        \foreach \from/\to in {1/2}
  \draw (\from) -- (\to);  
  \node (n5)[draw=white,fill=white,left of=1,node distance=1cm] {$\D$}; 
   \end{tikzpicture}
   \end{center}
 and assume we observe $\breve b_{31}> \breve b_{32} \breve b_{21}$. Then the matrix $\breve B$  fails to satisfy (\ref{eq:fixpoint}) and hence is not an element of $\calb(\D)$.  \halmos
\eexam
However, if we only estimate the ML coefficients corresponding to edges in $\D$ and then compute an estimate based on Lemma~\ref{lem:solving} below this phenomenon cannot occur.
 \ble\label{lem:solving}Let $B_0\in \overline{\mathbb{R}}_+^{d\times d}$ be a matrix with $b_{ji} >0\iff j\to i$. A matrix $A\in \overline{\mathbb{R}}_+^{d\times d}$ satisfies  
 \begin{equation}\label{eq:esteq}[a_{ji}>0 \iff j\in \An(i)] \text { and }A=I_d\vee (A\odot B_0)
\end{equation}
if and only if
 $A=(I_d\vee B_0)^{\odot (d-1)}$. 
 \ele
 \bproof 
We first show that $A=(I_d\vee B_0)^{\odot (d-1)}$ satisfies (\ref{eq:esteq}).   It is immediate that $a_{ji}>0 \iff j\in \An(i)$. We have (\citep{But2010}, Proposition 1.6.10) that
 \[(I_d \vee B_0)^{\odot (d-1)}=\bigvee_{k=0}^{d-1} B_0^{\odot k}=\bigvee_{k=0}^{\infty} B_0^{\odot k}\]
 and hence
 \[I_d\vee (A \odot B_0)=I_d\vee \{(I_d \vee B_0)^{\odot (d-1)}\odot B_0\}=I_d\vee \{\bigvee_{k=1}^{\infty} B_0^{\odot k}\}= \bigvee_{k=0}^{\infty} B_0^{\odot k}=A.\] 

It is easy to see directly that  $B_0^{\odot k}= 0$ for $k\geq d$ and hence if $\check A$ is a solution to \eqref{eq:esteq} we get by iteration, using that $(M\vee N)\odot K=(M\odot K)\vee (N\odot K)$,
\begin{eqnarray*}\check A &=& I_d \vee (\check A\odot B_0)\\& = &
I_d \vee [\{I_d\vee(\check A\odot B_0)\}\odot B_0]\\& = &I_d\vee B_0 \vee (\check A\odot B_0^{\odot 2})\\& =&\cdots\\&=& (I_d \vee B_0)^{\odot(d-1)} \vee (\check A \odot B_0^{\odot d}) =(I_d \vee B_0)^{\odot(d-1)}= A\end{eqnarray*}
and hence the solution to the equation is unique.
 \eproof

Thus we may define the estimate $\wh B$ by first calculating the matrix $\breve B_0= (\breve b_{ij} \bone_{\pa(j)}(i))_{d\times d}$ and then iterating the $\odot$-matrix product as: 
	\begin{align}\label{ch4:eq:Bhat}
\wh	B= (I_d\vee \breve B_0)^{\odot (d-1)}.
	\end{align}
It then follows that $\wh B_0=\breve B_0$  
	and  Lemma~\ref{lem:solving} yields that $\wh B$ is the unique element of $\mathcal{B}(\D)$ satisfying (\ref{ch4:eq:Bhat}).
By  Lemma~\ref{lem:XioverXj}(b), we also have 
\[b_{ji} \le\wh b_{ji} \le \breve b_{ji} \text{ for } j\in \an(i). \]
Consequently, when using $\wh B$ or $\breve B$ as an estimate of $B$, we never underestimate a \ML\ coefficient; furthermore, the matrix $\wh B$ always estimates $B$ more precisely than $\breve B$ and since we always have $\wh B\in \calb(\D)$, $\wh B$ seems to be clearly preferable as an estimate of $B$.

 \medskip
 
The following example shows how effective the estimate $\wh B$ can be; in particular,    $n$ does not necessarily need to be large.  
\bexam[One observation may be enough to estimate $B$ exactly\label{ch4:examBhat}]\\
 Consider the \DAG\ 
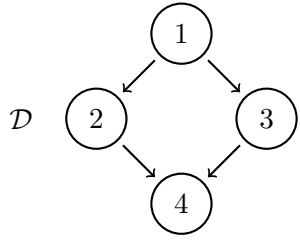
 \begin{center}
 \begin{tikzpicture}[->,every node/.style={circle,draw},line width=0.8pt, node distance=1.6cm,minimum size=0.8cm,outer sep=1mm]
  \node (1) [outer sep=1mm] {$1$};
   \node (2) [below left of=1] {$2$};
     \foreach \from/\to in {1/2}
   \draw (\from) -- (\to);   
   \node (3) [below right of=1] {$3$};
     \foreach \from/\to in {1/3}
   \draw (\from) -- (\to);   
   \node (4) [below right of=2] {$4$}; 
     \foreach \from/\to in {2/4,3/4}
   \draw (\from) -- (\to);   
    \node (n5)[draw=white,fill=white,left of=2,node distance=1cm] {$\D$}; 
 \end{tikzpicture}
 \end{center}
and assume that the paths $[1\to 2\to 4]$ and  $[1\to 3\to 4]$ are both max-weighted, which is equivalent to $b_{12}b_{24} = b_{13}b_{34}$.
If we observe the event 
\begin{align*} 
\big\{ X_2=b_{12}X_1\big\} \cap   \big\{ X_3=b_{13}X_1\big\} \cap  \big\{ X_4=b_{24}X_2\big\}  \cap \big\{ X_4=b_{34}X_3\big\}, 
\end{align*} 
then $\wh B=B$ so we estimate all \ML\ coefficients exactly. Note that this event has positive probability and occurs $\P$-almost surely  if and only if $Z_1$ realizes all node variables; i.e., if  $X_2=b_{12}Z_1$, $X_3=b_{13} Z_1$, and $X_4=b_{14}Z_1$. \halmos
\eexam

Since  by Table~\ref{ch4:table:rel} $\P(X_i=b_{ki}X_k)>0$ for $k\in\pa(i)$, it follows from the Borel-Cantelli lemma that $\wh b_{ki}$ $\P$-almost surely equals the true value for $n$ sufficiently large. Thus, if $n$ is large, $\wh B$ finds, with probability 1,  the true $B$.  In \cite{Davis1989} this is  discussed in a time-series framework used there and in \citet{davis:mccormick:89} they show that under suitable assumptions in the time-series framework, this estimator is asymptotically Fr\'echet distributed.
Assuming the probability of $\{X_i=b_{ki}X_k\}$ is known, we show next how one has to choose $n$ to observe this event with  probability greater than $1-p$ for some $p\in (0,1)$. We also prove that the probability for estimating the true $b_{ki}$ converges exponentially fast to 1.

\bpr \label{ch4:PropBhat}
Let   $\bfx^{(t)}=\big(X_1^{(t)},\ldots,X_n^{(t)}\big)$ for $t=1,\ldots,n$ be a sample from a recursive \ML\ model on a \DAG\ $\D$ with \ML\ coefficient matrix $B$. Let $i\in V$ and $k\in\pa(i)$.  
	 It then holds that 
	 \begin{align*}
	 \P\left(\bigwedge_{t=1}^nY_{ki}^{(t)}= b_{ki} \right)\ge 1-p\text{  for some $p\in(0,1)$ }
	 \end{align*}
	 if and only if  
	  \begin{align*}
	 n \ge   
	 \frac{\ln(p)}{\ln( \P(Y_{ki}> b_{ki}   ))}.
	 \end{align*}
Furthermore, 
the convergence  $\P\big(\bigwedge_{t=1}^n Y_{ki}^{(t)}= b_{ki} \big)\to 1$ 
 as $n\to \infty$ is exponentially fast. 
\epr
\bproof 
First note that the events $\{X_i=b_{ki}X_k\}$ and $\{X_i>b_{ki}X_k\}$ are complementary and both have positive probability. 
Further,
 using that $\bfx^{(1)},\ldots,\bfx^{(n)}$ are independent and identically distributed yields
	\begin{align*}
	\P\Big(\bigwedge_{t=1}^n Y_{ki}^{(t)}= b_{ki} \Big) = 1- \P\Big(\bigwedge_{t=1}^n Y_{ki}^{(t)}> b_{ki} \Big)
	=1-\prod_{t=1}^n\P({Y_{ki}^{(t)}}>b_{ki})=1- \P({Y_{ki}} > b_{ki})^n. 
	\end{align*}
Altogether, the statements follow.
\eproof
\noindent In conclusion,  $\wh B$ has the nice property to be  'geometrically consistent' in the sense that the probability of $\{\wh B=B\}$ converges exponentially fast to one.

\subsection*{The matrix $\boldsymbol{\wh B}$ is a generalized maximum likelihood estimate\label{ch4:s51}}
As we found in the previous section, the estimate $\wh B$ is preferable to the direct estimate $\breve B$ as it will always be closer to the true value. In this section we further establish that $\wh B$ is not just an \emph{ad hoc} estimator, but can indeed be derived from likelihood considerations.

For $B\in \calb(\D)$  and a fixed distribution of the innovation vector we let $P_B$ denote  the probability measure induced by a recursive \ML\ model  on $\D$ with \ML\ coefficient matrix $B$, i.e. the distribution of $\bfx$ where $\bfx=\bfz\odot B$. We shall denote the family of these probability measures by $\calp(\D)$.

We cannot use standard maximum likelihood methods to estimate  $B$, since the family $\calp(\D)$
is not dominated 
(cf. Example~4.4.1 of \cite{Gissibl2018}) 
and hence the standard likelihood function is not well defined. 
However, there exist generalizations of maximum likelihood estimation (GMLE) that cover the undominated case as well; \citet{Kalbfleisch1980},  Kiefer and Wolfowitz \citep{kiefer1956}, and \citet{Scholz1980} suggested such extensions.  
We essentially follow the Kiefer--Wolfowitz definition of a GMLE as also done, for example, by \citet{Gill1989} and \citet{Johansen1978}. 
In the following we shall show that $\wh B$ can be seen as a maximum likelihood estimate of $B$ in
 the extended sense introduced by Kiefer and Wolfowitz in \citep{kiefer1956}. 

Let $\calp$ be a family of probability measures on  $(\R_+^d, \borel(\R_+^d))$ where $\borel(\R_+^d)$ denotes the Borel $\sigma$-algebra on $\R_+^d$, 
and $\bfsx^{(1)},\ldots, \bfsx^{(n)}$ a random sample from some $P_0\in\calp$. For $P,Q\in\calp$ and $\bfsx\in\R_+^d$ we define \[\rho(\bfsx,P,Q):=\frac{dP}{d(P+Q)}(\bfsx),\] where ${dP}/{d(P+Q)}$ denotes a density of $P$ with respect to  $P+Q$. Then we call $\wh P$  a {\em generalized maximum likelihood estimate\/} of $P_0$ if
\begin{align}\label{ch4:def:gmle}
\prod_{t=1}^n \rho (\bfsx^{(t)},\widehat P, \widehat P)\neq 0 \quad \text{and} \quad  \prod_{t=1}^n\rho  (\bfsx^{(t)},Q, \widehat P) \le   
\prod_{t=1}^n \rho (\bfsx^{(t)},\widehat P, Q) \text{ for all $Q\in\calp$}.
\end{align}

Since $P$ is absolutely continuous with respect to $P+Q$, the density ${dP}/{d(P+Q)}$ always exists according to the Radon-Nikodym theorem. This means that the GMLE is well-defined, save for the usual ambiguity in the method of maximum likelihood that densities are only defined up to null sets and therefore a specific choice of densities must be made. 
The Kiefer--Wolfowitz definition extends the definition of a MLE in a very natural way as it simply says that for any $Q\in \calp$, $\wh P$ is the MLE in the smaller family $\{\hat P, Q\}$, consisting of only two distributions.
In \citep{kiefer1956}
only the second condition in \eqref{ch4:def:gmle} is required, but the first condition is implicit.  
The first step in verifying that $\wh B$ is a GMLE of $B$ is  to 
 specify densities of $P_B$ with respect to  $P_{B}+P_{B^*}$  for any two  $B, B^*\in\calb(\D)$. 
 For this purpose we determine a partition $\big\{A_0(B, B^*),  A_{1/2}(B, B^*)$, $A_1(B, B^*)\big\}$ of $\R_+^d$ that  satisfies the following  three properties,
 \begin{align}
	\text{\rm(A): }\quad& P_{ B}(A_0(B, B^*))=0,\nonumber\\
	\text{\rm(B): }\quad &P_{B} (A\cap A_{1/2}(B, B^*))=P_{B^*} (A\cap A_{1/2}(B, B^*))\text{ for every }A\in\borel(\R_+^d),\label{eq:partition}\\
	\text{ \rm (C): }\quad & P_{B^*}(A_1(B, B^*))=0.\nonumber
 \end{align}
Then we choose as density the measurable function from $\R_+^d$ to $\{0,1/2,1 \}$ defined as 
\begin{align}\label{ch4:dens:gen}
\bfsx\mapsto\rho(\bfsx, B, B^*):= 
\frac{1}{2}\cdot\bone_{A_{1/2}(B,B^*)}(\bfsx) + \bone_{A_1(B,B^*)}(\bfsx)=
\begin{cases}
0,  & \text{if $\bfsx\in A_0(B, B^*)$}, \\
\frac{1}{2},  & \text{if $\bfsx\in A_{1/2}(B, B^*)$}, \\
1,  & \text{if $\bfsx\in A_1(B, B^*)$}.
\end{cases}
\end{align}
This is a valid density because, using the properties (A), (B), (C), we obtain for every $A\in\borel(\R_+^d)$,
\begin{align*}
\int_A \rho(\bfsx, B,B^*) (P_{B}+P_{B^*})(d\bfsx)=  P_B(A\cap A_{1/2}(B,B^*))+  P_B(A\cap A_{1}(B,B^*))=P_{B}(A)\Nadine{.}
\end{align*}

We begin with an example that shall help to get an idea and provide insights into the concepts  and arguments we shall use in the general case. 
It is 
deliberately very detailed and although it deals with a very special case, it illustrates the main issues also for the general case.

\begin{figure}[htb]
	\centering 
	{\setlength{\tabcolsep}{-.45em}
		\renewcommand{\arraystretch}{0.0}
		\begin{tabular}{ccc}  
			\begin{tikzpicture}[scale=0.55]
			\begin{scope}
			\clip[postaction={fill=blue!40}](0,0)--(0,5)--(3.25,5);
			\end{scope}
			\begin{scope}
			\clip[postaction={fill=red!40}](0,0)--(5,0)--(5,3.25);
			\end{scope}
			\begin{scope}
			\clip[postaction={fill=OliveGreen!55}](0,0)--(3.25,5)--(5,5)--(5,3.25);
			\end{scope}
			\draw [OliveGreen!55,line width=2.5pt] (0, 0) -- (5,3.25);
			\draw
			(3,1.5)   node[OliveGreen,minimum size=2.4em,rotate=35] {\begin{scriptsize}$x_2= b_{12}x_1$\end{scriptsize}};
			\draw [red!40,line width=2.5pt] (0,0) -- (3.25,5);
			\draw
			(2.4,2.9)   node[red!80,rotate=58] {\begin{scriptsize}$x_2=b^*_{12}x_1$\end{scriptsize}}; 
			\draw [white,line width=1.5pt] (5, 0) -- (5,5);     
			\draw [white,line width=1.5pt] (0, 5) -- (5,5); 
			
			\draw [white,line width=1pt] (0, -0.05) -- (5,-0.05);
			\draw [white,line width=1pt] (-0.05, 0) -- (-0.05,5);
			\coordinate (y) at (0,5);
			\coordinate (x) at (-0.5,0);
			\draw[<-] (y) node[left] {\begin{footnotesize}$x_2$\end{footnotesize}} -- (0,-0.5);
			\draw[->] (x)-- (5,0) node[below]{\begin{footnotesize}$x_1$\end{footnotesize}};

			\draw [color=red!40,line width=3pt] (0, -1.375) -- (1.5,-1.375);
			\draw [color=blue!40,line width=3pt] (3.6,-1.375) -- (5,-1.375);
			\draw [OliveGreen!55,line width=3pt] (1.55, -1.375) -- (3.42,-1.375);
			\draw [red!40, fill=red!40]  (3.5,-1.375) circle (4pt) node {};
			\draw [OliveGreen!55, fill=OliveGreen!55]  (1.5,-1.375) circle (4pt) node {};
			\coordinate (x) at (-0.5,-2);
			\draw[->] (x)-- (5.5,-2) node[below] {\begin{small}$\frac{x_2}{x_1}$\end{small}};
			\draw (3.5,-2.1) -- (3.5,-1.9) node[below=2pt] {\begin{scriptsize}
				$b^*_{12}$ \end{scriptsize}}; 
			\draw (1.5,-2.1) -- (1.5,-1.9) node[below=2pt] {\begin{scriptsize}
				$b_{12}$
				\end{scriptsize}};

			\end{tikzpicture}
			&
			
			\begin{tikzpicture}[scale=0.55]
			\begin{scope}
			\clip[postaction={fill=blue!40}](0,0)--(0,5)--(5,5);
			\end{scope}
			\begin{scope}
			\clip[postaction={fill=red!40}](0,0)--(5,0)--(5,5);
			\end{scope}
			\draw [blue!40,line width=2.5pt] (0, 0) -- (5,5);
			\draw
			(3.15,2.5)   node[blue!80,rotate=45] {\begin{scriptsize}$x_2=b_{12}x_1=b^*_{12}x_1$\end{scriptsize}};
			\draw [white,line width=1.5pt] (5, 0) -- (5,5);     
			\draw [white,line width=1.5pt] (0, 5) -- (5,5); 
			
			\draw [white,line width=1pt] (0, -0.05) -- (5,-0.05);
			\draw [white,line width=1pt] (-0.05, 0) -- (-0.05,5);
			\coordinate (y) at (0,5);
			\coordinate (x) at (-0.5,0);
			\draw[<-] (y) node[left] {\begin{footnotesize}$x_2$\end{footnotesize}} -- (0,-0.5);
			\draw[->] (x)-- (5,0) node[below]{\begin{footnotesize}$x_1$\end{footnotesize}};
			
			\coordinate (x) at (-0.5,-2);
			\draw[->] (x)-- (5.5,-2) node[below] {\begin{small}$\frac{x_2}{x_1}$\end{small}};  
			\draw [blue!40,line width=3pt] (2.6,-1.375) -- (5,-1.375);
			\draw [red!40,line width=3pt] (0,-1.375) -- (2.5,-1.375);
			\draw [blue!40, fill=blue!40]  (2.5,-1.375) circle (4pt) node {};
			\draw (2.5,-2.1) -- (2.5,-1.9) node[below=2pt] {\begin{scriptsize}
				$b_{12}=b^*_{12}$ \end{scriptsize}};

			\end{tikzpicture}
			&
			\begin{tikzpicture}[scale=0.55]
			\begin{scope}
			\clip[postaction={fill=blue!40}](0,0)--(0,5)--(5,5);
			\end{scope}
			
			\begin{scope}
			\clip[postaction={fill=red!40}](0,0)--(5,0)--(5,5);
			\end{scope}
			\draw
			(3.15,2.5)   node[OliveGreen,rotate=45]  {\begin{scriptsize}$x_2=b_{12}x_1$\end{scriptsize}};
			\draw [red!60,line width=0.5pt] (0, 0) -- (5,2);
			\draw [line width=2.5pt,OliveGreen!60] (0, 0) -- (5,5);
			\draw
			(3.25,0.9)   node[red!80,rotate=25] {\begin{scriptsize}$x_2=b^*_{12}x_1$\end{scriptsize}};
			\draw [white,line width=1.5pt] (5, 0) -- (5,5);     
			\draw [white,line width=1.5pt] (0, 5) -- (5,5); 
			
			\draw [white,line width=1pt] (0, -0.05) -- (5,-0.05);
			\draw [white,line width=1pt] (-0.05, 0) -- (-0.05,5);
			\coordinate (y) at (0,5);
			\coordinate (x) at (-0.5,0);
			\draw[<-] (y) node[left] {\begin{footnotesize}$x_2$\end{footnotesize}} -- (0,-0.5);
			\draw[->] (x)-- (5,0) node[below]{\begin{footnotesize}$x_1$\end{footnotesize}};

			\coordinate (x) at (-0.5,-2);
			\draw[->] (x)-- (5.5,-2) node[below] {\begin{small}$\frac{x_2}{x_1}$\end{small}};   
			\draw [red!40,line width=3pt] (0, -1.375) -- (3.5,-1.375);
			\draw [blue!40,line width=3pt] (3.6, -1.375) -- (5,-1.375);
			\draw [OliveGreen!60, fill=OliveGreen!60]  (3.5,-1.375) circle (4pt) node {};
			\draw  (3.5,-2.1) -- (3.5,-1.9) node[below=2pt] {\begin{scriptsize}
				$b_{12}$ \end{scriptsize}}; 
			\draw  (1.5,-2.1) -- (1.5,-1.9) node[below=2pt] {\begin{scriptsize}
				$b^*_{12}$
				\end{scriptsize}};
			
			\end{tikzpicture}
	\end{tabular}}
	\caption{The density $\rho(\cdot,B,B^*)$ from Example~\ref{ch4:ex1} shown as a contour plot (top line) and as a function of $y_{12}={x_2}/{x_1}$ (bottom line) for the three situations $b_{12}<b^*_{12}$ (left-hand side), $b_{12}=b^*_{12}$ (middle), and $b_{12}>b^*_{12}$ (right-hand side). The area where it is ${\color{red!40}\boldsymbol{0}}/{\color{blue!40}\boldsymbol{\frac{1}{2}}}/{\color{OliveGreen!55}\boldsymbol{1}}$ is coloured in {\color{red!40}\bf red}/{\color{blue!40} \bf blue}/{\color{OliveGreen!55} \bf green}.\label{ch4:fig1:DAG12}}
\end{figure}
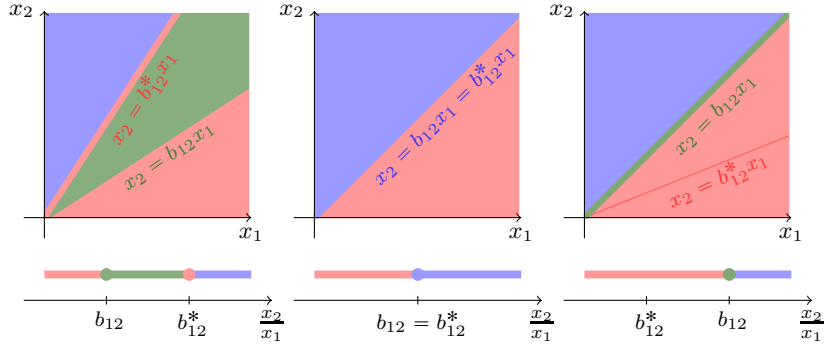

\bexam[How to find a density and the associated GMLEs\label{ch4:ex1}]\\
For $B, B^*\in\calb(\D)$ where $\D=(\{1,2\}, 1\to 2)$, we show that the partition 
\begin{align*}
\Big \{A_0(B, B^*)
&:=
\big\{\bfsx\in\R_+^2:  x_2 <b_{12} x_1 \big\}\cup \big\{\bfsx\in\R_+^2:    x_2=b^*_{12}x_1>b_{12}x_1\big\},\\
A_{1/2}(B,B^*) 
&:=
\big\{ \bfsx\in\R_+^2:  x_2=  b_{12}x_1=b^*_{12}x_1\big\}\cup  \big\{ \bfsx\in\R_+^2: x_2>    ( b_{12}\vee b^*_{12})x_1   \big\}, \\
A_1(B, B^*)&:= 
\big\{\bfsx\in\R_+^2:  b^*_{12}x_1 > x_2\ge  b_{12}x_1  \big\} \cup \big\{\bfsx\in\R_+^2:     x_2=b_{12}x_1 >b^*_{12}x_1 \big\}\Big \}
\end{align*}
of $\R_+^2$ satisfies properties (A), (B), (C) of (\ref{eq:partition}).   
Figure~\ref{ch4:fig1:DAG12} shows the corresponding density $\rho(\cdot,B,B^*)$ from \eqref{ch4:dens:gen}  for the three possible order relations between $b_{12}$ and $b^*_{12}$. 

Since by Table~\ref{ch4:table:rel},  $\ubpp({X_2}/{X_1})=[b_{12},\infty)$ and $b_{12}$ is the only atom of ${X_2}/{X_1}$, property (A) is true. By reversing  the roles of $B$ and $B^*$, (C) follows from (A).  
The condition (B) is obvious if $b_{12}=b^*_{12}$. 
Assume that $b_{12}\neq b_{12}^*$.  We then have by definition of $\bfx$ that  $\{\bfx\in A_{1/2}(B, B^*)  \}=\{  X_2>(b_{12}\vee b^*_{12})  X_1 \}=\{  Z_2> (b_{12}\vee b^*_{12}) Z_1  \}$ and  $ X_2=Z_2$ on $\{ Z_2> (b_{12}\vee b^*_{12}) Z_1  \}$. With this, using  that 
$A_{1/2}(B^*, B)=A_{1/2}(B, B^*) $,
we obtain for  $A\in\mathbb{B}(\R_+^2)$, 
\begin{align*}
P_{B} (A\cap A_{1/2}(B, B^*))&=\P(\{ \bfx\in A\} \cap \{Z_2> (b_{12}\vee b^*_{12})  Z_1\}) \\
&=\P(\{  ( Z_1,  Z_2) \in A\} \cap \{ Z_2> (b_{12}\vee b^*_{12})  Z_1\} )\\
&= P_{B^*}(A\cap A_{1/2}(B^*,B))=P_{B^*}(A\cap A_{1/2}(B,B^*)).
\end{align*} 
 
We now use the density found to determine the GMLE of $B$. 
The only \ML\ coefficient we have to estimate is $b_{12}$.  
As before we let $\wh b_{12}=\breve b_{12}$ be the minimal observed ratio of ${X_2}/{X_1}$ 
and let $\wh B$ be the corresponding \ML\ coefficient matrix from (\ref{ch4:eq:Bhat}). 
Defining $n(B,B^*)=\vert \{t: \bfsx^{(t)}\in A_{1/2}(B,B^*) \}\vert$ and using that $n(B,B^*)=n(B^*,B)$, we obtain
\begin{align*}
\prod_{t=1}^n \rho(\bfsx^{(t)},B,B^*)&=2^{-n(B,B^*)} \prod_{t=1}^n \bone_{\R_+^d\setminus A_0(B,B^*)}\big(\bfsx^{(t)}\big),\\
\prod_{t=1}^n \rho(\bfsx^{(t)},B^*,B)&=2^{-n(B,B^*)} \prod_{t=1}^n \bone_{\R_+^d\setminus A_0(B^*,B)}\big(\bfsx^{(t)}\big). 
\end{align*}
Let now $\wt B$ be an arbitrary potential GMLE of $B$. Then $P_{\wt B}\in\calp(\D)$ satisfies the first condition in \eqref{ch4:def:gmle} if and only if 
\begin{align}\label{ch4:cond1:DAG12}
\wt b_{12} x_1^{(t)}\le x_2^{(t)} \text{ for all $t$, equivalently $\wt b_{12}\le \wh b_{12}$}
\end{align} 
and the second condition if and only if
\begin{align}\label{ch4:cond2:DAG12}
\text{for all $B\in\calb(\D)$, if  some $\bfsx^{(t)}\in A_0(\wt B,B)$, then some $\bfsx^{(s)}\in A_0(B, \wt B)$.}
\end{align}
In summary, some $\wt B\in \mathcal{B}(\D)$ is a GMLE of $ B$ if and only if \eqref{ch4:cond1:DAG12} and  \eqref{ch4:cond2:DAG12} are satisfied. We discuss the possible GMLEs of $b_{12}$ in detail. 
\begin{enumerate}
	\item[(a)] $\wt b_{12}<\wh b_{12}$ is no GMLE:\\
	Set $b_{12}=\wh b_{12}$, and let $\bfsx^{(t)}$ be such that $\wh b_{12}x_1^{(t)}=x_2^{(t)}$. Then $\bfsx^{(t)}\in \big\{ \bfsx\in\R_+^2: x_2=b_{12}x_1 >\wt b_{12} x_2\big\}\subseteq A_0(\wt B, B)$ but no $\bfsx^{(s)} \in A_0(B,\wt B)= \big\{ \bfsx\in\R_+^2: x_2<b_{12}x_1 \big\}$. This contradicts \eqref{ch4:cond2:DAG12}; consequently, $\wt b_{12}$ cannot be a GMLE of $b_{12}$. In 
	Figure~\ref{ch4:fig2:DAG12}(a) we illustrate this situation. On the left-hand side a contour plot of the density $\rho(\cdot, \wt B, B)$ is shown, on the right-hand side of $\rho(\cdot,  B,\wt B)$.  The crosses represent the realizations $\bfsx^{(1)},\ldots, \bfsx^{(n)}$. In the left plot crosses are in the $0$-area coloured in red, namely, those that realize $\wh b_{12}$, but in the right plot not.  So $\wt B$ cannot be a GMLE of $B$. 
	\item[(b)] $\wt b_{12}>\wh b_{12}$ is no GMLE:\\ This follows directly from (\ref{ch4:cond1:DAG12}).
Figure~\ref{ch4:fig2:DAG12}(b) shows a situation that contradicts \eqref{ch4:cond2:DAG12}, similarly to Figure~\ref{ch4:fig2:DAG12}(a) in (1). 
	\item[(c)] $\wt b_{12}=\wh b_{12}$ is a GMLE:\\
	Condition \eqref{ch4:cond1:DAG12}  holds obviously. To prove \eqref{ch4:cond2:DAG12}, assume for some $B\in\calb(\D)$ that some $\bfsx^{(t)}\in A_0(\wh B, B)$. By definition of $A_0(\wh B,B)$, $x_2^{(t)}=b_{12}x_1^{(t)}> \wh b_{12} x_1^{(t)}$, which implies that $b_{12}>\wh b_{12}$.  For $\bfsx^{(s)}$ such that $\wh b_{12}x_1^{(s)}=x_2^{(s)}$, we then find that $ x_2^{(s)}< b_{12}x_1^{(s)}$. Hence, $\bfsx^{(s)}\in A_0(B,\wh B)$, and $\wh b_{12}$ is a GMLE of $b_{12}$. We learn this informally from Figure~\ref{ch4:fig2:DAG12}(c). The top line shows contour plots  of $\rho(\cdot, \wh B,B)$ for the three different orders between $b_{12}$ and $\wh b_{12}$, and the bottom line shows the corresponding contour plots of  $\rho(\cdot,  B, \wh B)$. The two plots on the left-hand side  correspond to  the situation from above: in the upper plot there are realizations in the $0$-area, namely those that are on the line $\big\{ \bfsx\in\R_+^2: x_2=b_{12}x_1\big\}$, but then there are also realizations in the $0$-area of the lower plot (those that lie below this line).  Hence, \eqref{ch4:cond2:DAG12} holds. Since there is no realization in the $0$-area of the middle and right plot in the top line,  \eqref{ch4:cond2:DAG12}  is automatically satisfied if $b_{12}\le \wh b_{12}$. 
	\begin{figure}
		\centering
		\subfloat[$\wt b_{12}<\wh b_{12}$ is no  GMLE.]{
			{\setlength{\tabcolsep}{-.45em}
				\renewcommand{\arraystretch}{0.0}
				\begin{tabular}{cc}
					\begin{tikzpicture}[scale=0.55]  
					\begin{scope}
					\clip[postaction={fill=blue!40}](0,0)--(0,5)--(5,5);
					\end{scope}
					\begin{scope}
					\clip[postaction={fill=red!40}](0,0)--(5,2) -- (5,0);
					\end{scope}
					\begin{scope}
					\clip[postaction={fill=OliveGreen!55}](0,0)--(5,5)--(5,2) ;
					\end{scope}
					\draw [red!40,line width=2.5pt] (0, 0) -- (5,5);
					\draw  (3.15,2.5)   node[red!80,rotate=45,minimum size=2.4em] {\begin{scriptsize}$x_2=  b_{12}x_1=\wh b_{12}x_1$\end{scriptsize}};
					\draw (1.5,1.5) node[cross=1.75pt] {};
					\draw (1.5,3.333) node[cross=1.75pt] {};
					\draw (0.5,1.111) node[cross=1.75pt] {};
					\draw (4,4) node[cross=1.75pt] {};
					\draw (3,3) node[cross=1.75pt] {};
					\draw (1,4) node[cross=1.75pt] {};
					\draw (2.5,4) node[cross=1.75pt] {};
					\draw (0.5,3) node[cross=1.75pt] {};
					\draw (2,3.5) node[cross=1.75pt] {};
					\draw (0.5,4.5) node[cross=1.75pt] {};
					\draw (1.75,2.5) node[cross=1.75pt] {};
					\draw (3.3,0.8)   node[OliveGreen,rotate=25,minimum size=2.4em] {\begin{scriptsize}$x_2= \wt b_{12}x_1$\end{scriptsize}};
					\draw [OliveGreen!55,line width=2.5pt] (0, 0) -- (5,2);
					\draw (2.5,5)   node[above=1pt] {\begin{scriptsize}
						$\rho(\cdot,\wt B,B)$
						\end{scriptsize}};
					\draw [white,line width=1.5pt] (5, 0) -- (5,5);     
					\draw [white,line width=1.5pt] (0, 5) -- (5,5); 
					
					\draw [white,line width=1pt] (0, -0.05) -- (5,-0.05);
					\draw [white,line width=1pt] (-0.05, 0) -- (-0.05,5);
					\coordinate (y) at (0,5);
					\coordinate (x) at (-0.5,0);
					\draw[<-] (y) node[left] {\begin{footnotesize}$x_2$\end{footnotesize}} -- (0,-0.5);
					\draw[->] (x)-- (5,0) node[below]{\begin{footnotesize}$x_1$\end{footnotesize}};  
					\end{tikzpicture}

					&
					\begin{tikzpicture}[scale=0.55]  
					\begin{scope}
					\clip[postaction={fill=blue!40}](0,0)--(0,5)--(5,5);
					\end{scope}
					\begin{scope}
					\clip[postaction={fill=red!40}](0,0)--(5,0)--(5,5);
					\end{scope}
					\draw [OliveGreen!60,line width=2.5pt] (0, 0) -- (5,5);
					\draw (1.5,1.5) node[cross=1.75pt] {};
					\draw (1.5,3.333) node[cross=1.75pt] {};
					\draw (0.5,1.111) node[cross=1.75pt] {};
					\draw (4,4) node[cross=1.75pt] {};
					\draw (3,3) node[cross=1.75pt] {};
					\draw (1,4) node[cross=1.75pt] {};
					\draw (2.5,4) node[cross=1.75pt] {};
					\draw (0.5,3) node[cross=1.75pt] {};
					\draw (2,3.5) node[cross=1.75pt] {};
					\draw (0.5,4.5) node[cross=1.75pt] {};
					\draw (1.75,2.5) node[cross=1.75pt] {};
					\draw  (3.15,2.5)   node[OliveGreen,rotate=45,minimum size=2.4em] {\begin{scriptsize}$x_2=  b_{12}x_1=\wh b_{12}x_1$\end{scriptsize}};
					\draw [red!60,line width=0.5pt] (0, 0) -- (5,2); 
					\draw (2.5,5)  node [above=1pt]  {\begin{scriptsize}
						$\rho(\cdot, B, \wt B)$
						\end{scriptsize}};
					\draw [white,line width=1.5pt] (5, 0) -- (5,5);     
					\draw [white,line width=1.5pt] (0, 5) -- (5,5); 
					
					\draw [white,line width=1pt] (0, -0.05) -- (5,-0.05);
					\draw [white,line width=1pt] (-0.05, 0) -- (-0.05,5);
					\coordinate (y) at (0,5);
					\coordinate (x) at (-0.5,0);
					\draw[<-] (y) node[left] {\begin{footnotesize}$x_2$\end{footnotesize}} -- (0,-0.5);
					\draw[->] (x)-- (5,0) node[below]{\begin{footnotesize}$x_1$\end{footnotesize}}; 
					\draw (3.25,0.9)  node[red!80,rotate=25,minimum size=2.4em] {\begin{scriptsize}$x_2= \wt b_{12}x_1$\end{scriptsize}}; 
					\end{tikzpicture}
				\end{tabular}
		}  }
		\subfloat[$\wt b_{12}>\wh b_{12}$ is no GMLE.]{
			{\setlength{\tabcolsep}{-.45em}
				\renewcommand{\arraystretch}{0.0}
				\begin{tabular}{cc}
					\begin{tikzpicture}[scale=0.55]
					\begin{scope}
					\clip[postaction={fill=blue!40}](0,0)--(0,5)--(4.1,5);
					\end{scope}
					\begin{scope}
					\clip[postaction={fill=red!40}](0,0)--(5,0)--(5,5) -- (4.1,5);
					\end{scope}
					\draw  (2.65,2.5) node[OliveGreen,minimum size=2.4em,rotate=55] {\begin{scriptsize}$x_2= \wt b_{12}x_1$\end{scriptsize}};  
					\draw [red!60,line width=0.5pt] (0, 0) -- (5,2);
					\draw  (3.25,0.9)   node[red!80,rotate=25] {\begin{scriptsize}$x_2=b_{12}x_1$\end{scriptsize}};
					\draw [line width=2.5pt,OliveGreen!60] (0, 0) -- (4.1,5);
					\draw (1.5,1.5) node[cross=1.75pt] {};
					\draw (1.5,3.333) node[cross=1.75pt] {};
					\draw (0.5,1.111) node[cross=1.75pt] {};
					\draw (4,4) node[cross=1.75pt] {};
					\draw (3,3) node[cross=1.75pt] {};
					\draw (1,4) node[cross=1.75pt] {};
					\draw (2.5,4) node[cross=1.75pt] {};
					\draw (0.5,3) node[cross=1.75pt] {};
					\draw (2,3.5) node[cross=1.75pt] {};
					\draw (0.5,4.5) node[cross=1.75pt] {};
					\draw (1.75,2.5) node[cross=1.75pt] {};
					\draw [white,line width=1.5pt] (5, 0) -- (5,5);     
					\draw [white,line width=1.5pt] (0, 5) -- (5,5); 
					
					\draw [white,line width=1pt] (0, -0.05) -- (5,-0.05);
					\draw [white,line width=1pt] (-0.05, 0) -- (-0.05,5);
					\coordinate (y) at (0,5);
					\coordinate (x) at (-0.5,0);
					\draw[<-] (y) node[left] {\begin{footnotesize}$x_2$\end{footnotesize}} -- (0,-0.5);
					\draw[->] (x)-- (5,0) node[below]{\begin{footnotesize}$x_1$\end{footnotesize}}; 
					\draw (2.5,5)   node[above=1pt] {\begin{scriptsize}
						$\rho(\cdot,\wt B,B)$
						\end{scriptsize}}; 
					\end{tikzpicture}

					&
					\begin{tikzpicture}[scale=0.55]
					\begin{scope}
					\clip[postaction={fill=blue!40}](0,0)--(0,5)--(4.1,5);
					\end{scope}
					\begin{scope}
					\clip[postaction={fill=red!40}](0,0)--(5,0)--(5,5) -- (4.1,5);
					\end{scope}
					\begin{scope}
					\clip[postaction={fill=OliveGreen!55}](0,0)--(4.1,5) -- (5,5)--(5,2);
					\end{scope}
					\draw  (2.65,2.5)  node[red!80,minimum size=2.4em,rotate=55] {\begin{scriptsize}$x_2= \wt b_{12}x_1$\end{scriptsize}};
					\draw [OliveGreen!55,line width=2.5pt] (0, 0) -- (5,2);
					\draw [line width=2.5pt,red!40] (0, 0) -- (4.1,5);    
					\draw    (3.3,0.8)  node[OliveGreen,rotate=25] {\begin{scriptsize}$x_2=b_{12}x_1$\end{scriptsize}};
					\draw (1.5,1.5) node[cross=1.75pt] {};
					\draw (1.5,3.333) node[cross=1.75pt] {};
					\draw (0.5,1.111) node[cross=1.75pt] {};
					\draw (4,4) node[cross=1.75pt] {};
					\draw (3,3) node[cross=1.75pt] {};
					\draw (1,4) node[cross=1.75pt] {};
					\draw (2.5,4) node[cross=1.75pt] {};
					\draw (0.5,3) node[cross=1.75pt] {};
					\draw (2,3.5) node[cross=1.75pt] {};
					\draw (0.5,4.5) node[cross=1.75pt] {};
					\draw (1.75,2.5) node[cross=1.75pt] {};
					\draw [white,line width=1.5pt] (5, 0) -- (5,5);     
					\draw [white,line width=1.5pt] (0, 5) -- (5,5); 
					
					\draw [white,line width=1pt] (0, -0.05) -- (5,-0.05);
					\draw [white,line width=1pt] (-0.05, 0) -- (-0.05,5);
					\coordinate (y) at (0,5);
					\coordinate (x) at (-0.5,0);
					\draw[<-] (y) node[left] {\begin{footnotesize}$x_2$\end{footnotesize}} -- (0,-0.5);
					\draw[->] (x)-- (5,0) node[below]{\begin{footnotesize}$x_1$\end{footnotesize}};  
					\draw (2.5,5)   node[above=1pt] {\begin{scriptsize}
						$\rho(\cdot, B,\wt B)$
						\end{scriptsize}};
					\end{tikzpicture}
				\end{tabular}
		}  }\\
		\subfloat[$\wt b_{12}=\wh b_{12}$ is a GMLE.]{
			{ \hspace*{-5em}
				\setlength{\tabcolsep}{-.45em}
				\renewcommand{\arraystretch}{0.0}
				\begin{tabular}{ccc}
					\begin{tikzpicture}[scale=0.55]  
					\draw (-2,3) node {\begin{scriptsize}
						$\rho(\cdot,\wh B, B)$:
						\end{scriptsize}};
					\begin{scope}
					\clip[postaction={fill=blue!40}](0,0)--(0,5)--(2.25,5);
					\end{scope}
					\begin{scope}
					\clip[postaction={fill=red!40}](0,0)--(5,0)--(5,5);
					\end{scope}
					\begin{scope}
					\clip[postaction={fill=OliveGreen!55}](0,0)--(2.25,5)--(5,5);
					\end{scope}
					\draw [OliveGreen!55,line width=2.5pt] (0, 0) -- (5,5);
					\draw  (3.15,2.5)  node[OliveGreen,minimum size=2.4em,rotate=45] {\begin{scriptsize}$x_2= \wh b_{12}x_1$\end{scriptsize}};
					\draw [red!40,line width=2.5pt] (0,0) -- (2.25,5);
					\draw  (1,3.2)   node[red!80,rotate=67] {\begin{scriptsize}$x_2=b_{12}x_1$\end{scriptsize}}; 
					\draw (1.5,1.5) node[cross=1.75pt] {};
					\draw (1.5,3.333) node[cross=1.75pt] {};
					\draw (0.5,1.111) node[cross=1.75pt] {};
					\draw (4,4) node[cross=1.75pt] {};
					\draw (3,3) node[cross=1.75pt] {};
					\draw (1,4) node[cross=1.75pt] {};
					\draw (2.5,4) node[cross=1.75pt] {};
					\draw (0.5,3) node[cross=1.75pt] {};
					\draw (2,3.5) node[cross=1.75pt] {};
					\draw (0.5,4.5) node[cross=1.75pt] {};
					\draw (1.75,2.5) node[cross=1.75pt] {};
					\draw [white,line width=1.5pt] (5, 0) -- (5,5);     
					\draw [white,line width=1.5pt] (0, 5) -- (5,5); 
					
					\draw [white,line width=1pt] (0, -0.05) -- (5,-0.05);
					\draw [white,line width=1pt] (-0.05, 0) -- (-0.05,5);
					\coordinate (y) at (0,5);
					\coordinate (x) at (-0.5,0);
					\draw[<-] (y) node[left] {\begin{footnotesize}$x_2$\end{footnotesize}} -- (0,-0.5);
					\draw[->] (x)-- (5,0) node[below]{\begin{footnotesize}$x_1$\end{footnotesize}};  
					\end{tikzpicture}
					&
					\begin{tikzpicture}[scale=0.55]

					\begin{scope}
					\clip[postaction={fill=blue!40}](0,0)--(0,5)--(5,5);
					\end{scope}
					\begin{scope}
					\clip[postaction={fill=red!40}](0,0)--(5,0)--(5,5);
					\end{scope}
					\draw [blue!40,line width=2.5pt] (0, 0) -- (5,5);
					
					\draw (3.15,2.5)   node[blue!80,rotate=45] {\begin{scriptsize}$x_2=\wh b_{12}x_1=b_{12}x_1$\end{scriptsize}};
					\draw (1.5,1.5) node[cross=1.75pt] {};
					\draw (1.5,3.333) node[cross=1.75pt] {};
					\draw (0.5,1.111) node[cross=1.75pt] {};
					\draw (4,4) node[cross=1.75pt] {};
					\draw (3,3) node[cross=1.75pt] {};
					\draw (1,4) node[cross=1.75pt] {};
					\draw (2.5,4) node[cross=1.75pt] {};
					\draw (0.5,3) node[cross=1.75pt] {};
					\draw (2,3.5) node[cross=1.75pt] {};
					\draw (0.5,4.5) node[cross=1.75pt] {};
					\draw (1.75,2.5) node[cross=1.75pt] {};
					\draw [white,line width=1.5pt] (5, 0) -- (5,5);     
					\draw [white,line width=1.5pt] (0, 5) -- (5,5); 
					
					\draw [white,line width=1pt] (0, -0.05) -- (5,-0.05);
					\draw [white,line width=1pt] (-0.05, 0) -- (-0.05,5);
					\coordinate (y) at (0,5);
					\coordinate (x) at (-0.5,0);
					\draw[<-] (y) node[left] {\begin{footnotesize}$x_2$\end{footnotesize}} -- (0,-0.5);
					\draw[->] (x)-- (5,0) node[below]{\begin{footnotesize}$x_1$\end{footnotesize}};  
					\end{tikzpicture}
					&
					\begin{tikzpicture}[scale=0.55]
					\begin{scope}
					\clip[postaction={fill=blue!40}](0,0)--(0,5)--(5,5);
					\end{scope}
					\begin{scope}
					\clip[postaction={fill=red!40}](0,0)--(5,0)--(5,5);
					\end{scope}
					\draw (3.15,2.5)  node[OliveGreen,rotate=45,minimum size=2.4em] {\begin{scriptsize}$x_2=\wh b_{12}x_1$\end{scriptsize}};
					\draw [red!60,line width=0.5pt] (0, 0) -- (5,2);         
					\draw  (3.25,0.9)     node[red!80,rotate=25]  {\begin{scriptsize}$x_2=b_{12}x_1$\end{scriptsize}};
					\draw [line width=2.5pt,OliveGreen!60] (0, 0) -- (5,5);
					\draw (1.5,1.5) node[cross=1.75pt] {};
					\draw (1.5,3.333) node[cross=1.75pt] {};
					\draw (0.5,1.111) node[cross=1.75pt] {};
					\draw (4,4) node[cross=1.75pt] {};
					\draw (3,3) node[cross=1.75pt] {};
					\draw (1,4) node[cross=1.75pt] {};
					\draw (2.5,4) node[cross=1.75pt] {};
					\draw (0.5,3) node[cross=1.75pt] {};
					\draw (2,3.5) node[cross=1.75pt] {};
					\draw (0.5,4.5) node[cross=1.75pt] {};
					\draw (1.75,2.5) node[cross=1.75pt] {};
					\draw [white,line width=1.5pt] (5, 0) -- (5,5);     
					\draw [white,line width=1.5pt] (0, 5) -- (5,5); 
					
					\draw [white,line width=1pt] (0, -0.05) -- (5,-0.05);
					\draw [white,line width=1pt] (-0.05, 0) -- (-0.05,5);
					\coordinate (y) at (0,5);
					\coordinate (x) at (-0.5,0);
					\draw[<-] (y) node[left] {\begin{footnotesize}$x_2$\end{footnotesize}} -- (0,-0.5);
					\draw[->] (x)-- (5,0) node[below]{\begin{footnotesize}$x_1$\end{footnotesize}};  
					\end{tikzpicture}\\
					\begin{tikzpicture}[scale=0.55]  
					\draw (-2,3) node {\begin{scriptsize}
						$\rho(\cdot,B,\wh B)$:
						\end{scriptsize}};
					\begin{scope}
					\clip[postaction={fill=blue!40}](0,0)--(0,5)--(2.25,5);
					\end{scope}
					\begin{scope}
					\clip[postaction={fill=red!40}](0,0)--(5,0)--(5,5);
					\end{scope}
					\begin{scope}
					\clip[postaction={fill=red!40}](0,0)--(2.25,5)--(5,5);
					\end{scope}
					\draw [red!40,line width=2.5pt] (0, 0) -- (5,5);
					\draw   (3.15,2.5)  node[red!80,minimum size=2.4em,rotate=45] {\begin{scriptsize}$x_2= \wh b_{12}x_1$\end{scriptsize}};
					\draw [OliveGreen!60,line width=2.5pt] (0,0) -- (2.25,5);
					\draw [red!60,line width=0.5pt] (0,0) -- (5,5);
					\draw  (1,3.2)    node[OliveGreen,rotate=67] {\begin{scriptsize}$x_2=b_{12}x_1$\end{scriptsize}}; 
					\draw (1.5,1.5) node[cross=1.75pt] {};
					\draw (1.5,3.333) node[cross=1.75pt] {};
					\draw (0.5,1.111) node[cross=1.75pt] {};
					\draw (4,4) node[cross=1.75pt] {};
					\draw (3,3) node[cross=1.75pt] {};
					\draw (1,4) node[cross=1.75pt] {};
					\draw (2.5,4) node[cross=1.75pt] {};
					\draw (0.5,3) node[cross=1.75pt] {};
					\draw (2,3.5) node[cross=1.75pt] {};
					\draw (0.5,4.5) node[cross=1.75pt] {};
					\draw (1.75,2.5) node[cross=1.75pt] {};
					\draw [white,line width=1.5pt] (5, 0) -- (5,5);     
					\draw [white,line width=1.5pt] (0, 5) -- (5,5); 
					
					\draw [white,line width=1pt] (0, -0.05) -- (5,-0.05);
					\draw [white,line width=1pt] (-0.05, 0) -- (-0.05,5);
					\coordinate (y) at (0,5);
					\coordinate (x) at (-0.5,0);
					\draw[<-] (y) node[left] {\begin{footnotesize}$x_2$\end{footnotesize}} -- (0,-0.5);
					\draw[->] (x)-- (5,0) node[below]{\begin{footnotesize}$x_1$\end{footnotesize}};  
					\end{tikzpicture}
					&
					
					\begin{tikzpicture}[scale=0.55]
					\begin{scope}
					\clip[postaction={fill=blue!40}](0,0)--(0,5)--(5,5);
					\end{scope}
					\begin{scope}
					\clip[postaction={fill=red!40}](0,0)--(5,0)--(5,5);
					\end{scope}
					\draw [blue!40,line width=2.5pt] (0, 0) -- (5,5);
					\draw  (3.15,2.5)  node[blue!80,rotate=45] {\begin{scriptsize}$x_2=\wh b_{12}x_1=b_{12}x_1$\end{scriptsize}};
					\draw (1.5,1.5) node[cross=1.75pt] {};
					\draw (1.5,3.333) node[cross=1.75pt] {};
					\draw (0.5,1.111) node[cross=1.75pt] {};
					\draw (4,4) node[cross=1.75pt] {};
					\draw (3,3) node[cross=1.75pt] {};
					\draw (1,4) node[cross=1.75pt] {};
					\draw (2.5,4) node[cross=1.75pt] {};
					\draw (0.5,3) node[cross=1.75pt] {};
					\draw (2,3.5) node[cross=1.75pt] {};
					\draw (0.5,4.5) node[cross=1.75pt] {};
					\draw (1.75,2.5) node[cross=1.75pt] {};
					\draw [white,line width=1.5pt] (5, 0) -- (5,5);     
					\draw [white,line width=1.5pt] (0, 5) -- (5,5); 
					
					\draw [white,line width=1pt] (0, -0.05) -- (5,-0.05);
					\draw [white,line width=1pt] (-0.05, 0) -- (-0.05,5);
					\coordinate (y) at (0,5);
					\coordinate (x) at (-0.5,0);
					\draw[<-] (y) node[left] {\begin{footnotesize}$x_2$\end{footnotesize}} -- (0,-0.5);
					\draw[->] (x)-- (5,0) node[below]{\begin{footnotesize}$x_1$\end{footnotesize}};  
					\end{tikzpicture}
					&
					\begin{tikzpicture}[scale=0.55]
					\begin{scope}
					\clip[postaction={fill=blue!40}](0,0)--(0,5)--(5,5);
					\end{scope}
					\begin{scope}
					\clip[postaction={fill=red!40}](0,0)--(5,0)--(5,5);
					\end{scope}
					\begin{scope}
					\clip[postaction={fill=OliveGreen!55}](0,0)--(5,5)--(5,2);
					\end{scope}
					\draw [OliveGreen!55,line width=2.5pt] (0, 0) -- (5,2);
					\draw [line width=2.5pt,red!40] (0, 0) -- (5,5);
					\draw  (3.15,2.5)  node[red!80,rotate=45,minimum size=2.4em] {\begin{scriptsize}$x_2=\wh b_{12}x_1$\end{scriptsize}};
					\draw   (3.3,0.8)   node[OliveGreen,rotate=25] {\begin{scriptsize}$x_2= b_{12}x_1$\end{scriptsize}};
					\draw (1.5,1.5) node[cross=1.75pt] {};
					\draw (1.5,3.333) node[cross=1.75pt] {};
					\draw (0.5,1.111) node[cross=1.75pt] {};
					\draw (4,4) node[cross=1.75pt] {};
					\draw (3,3) node[cross=1.75pt] {};
					\draw (1,4) node[cross=1.75pt] {};
					\draw (2.5,4) node[cross=1.75pt] {};
					\draw (0.5,3) node[cross=1.75pt] {};
					\draw (2,3.5) node[cross=1.75pt] {};
					\draw (0.5,4.5) node[cross=1.75pt] {};
					\draw (1.75,2.5) node[cross=1.75pt] {};
					\draw [white,line width=1.5pt] (5, 0) -- (5,5);     
					\draw [white,line width=1.5pt] (0, 5) -- (5,5); 
					
					\draw [white,line width=1pt] (0, -0.05) -- (5,-0.05);
					\draw [white,line width=1pt] (-0.05, 0) -- (-0.05,5);
					\coordinate (y) at (0,5);
					\coordinate (x) at (-0.5,0);
					\draw[<-] (y) node[left] {\begin{footnotesize}$x_2$\end{footnotesize}} -- (0,-0.5);
					\draw[->] (x)-- (5,0) node[below]{\begin{footnotesize}$x_1$\end{footnotesize}};  
					\end{tikzpicture}
			\end{tabular}
		} }\\

		\caption{Discussion  of the  GMLEs of $b_{12}$ with respect to the density 
			from Figure~\ref{ch4:fig1:DAG12}.\label{ch4:fig2:DAG12}; see further explanation in (a), (b), and (c) of Example~\ref{ch4:ex1}. }
	\end{figure}
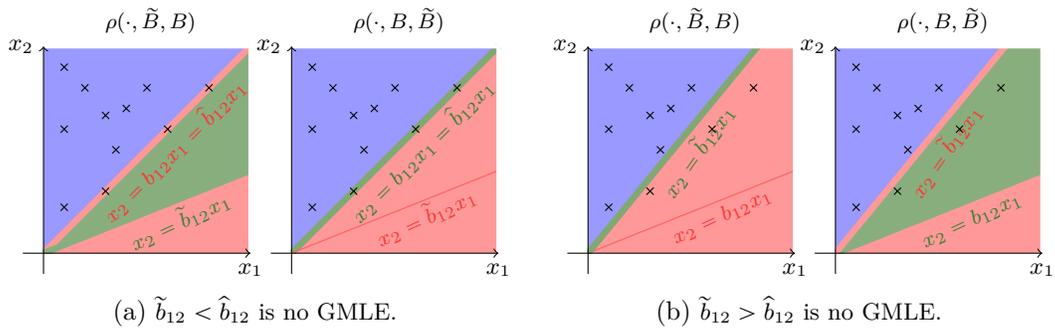
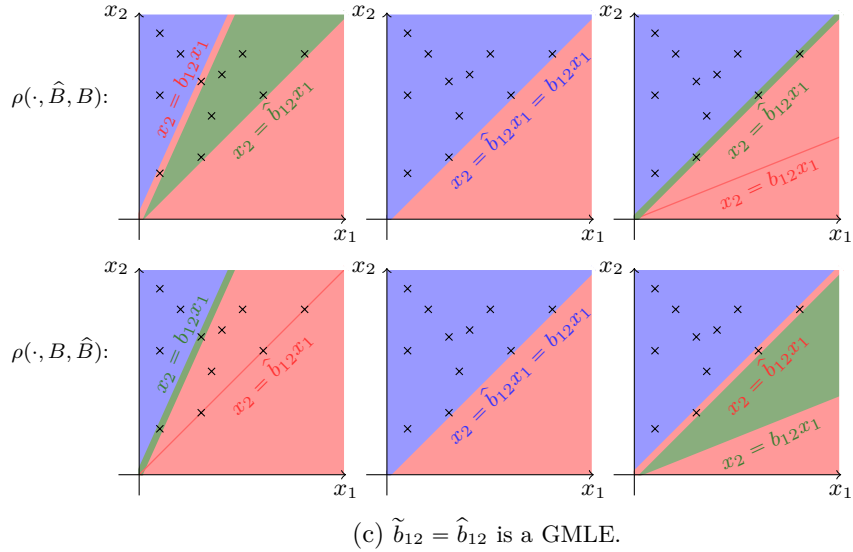
	\halmos
\end{enumerate}
\eexam

	In what follows we specify, for the general case,  one density of $P_B$ with respect to $P_B+P_{B^*}$  that has a representation as in \eqref{ch4:dens:gen} and leads to $\wh B$ as a GMLE of $B$.

	Our partition  $\big\{A_0(B,B^*), A_{1/2}(B,B^*), A_1(B,B^*)\big\}$ of $\R_+^d$  is based on the following representation for the components of $\bfx$:
	\begin{align}\label{ch4:proof:dens}
	X_i=\bigvee_{k\in\pa(i)} b_{ki}X_k\vee Z_i; \quad \text{in particular, $ X_i\ge \bigvee_{k\in\pa(i)}b_{ki}X_k$}, \quad i\in V.
	\end{align} 
	We begin with the specification of $A_{1/2}(B,B^*)$ and prove a property needed subsequently to verify property (B). Have in mind 
	that if $b_{ki}>b^*_{ki}$  for all $k\in\pa(i)$ or $b_{ki}<b^*_{ki}$ for all $k\in\pa(i)$ then $\big\{\bfsx\in\R_+^d: x_{i}  = \bigvee_{k\in \pa(i)} b^*_{ki}x_k=\bigvee_{k\in \pa(i)}  b_{k i}x_k  \big\}=\emptyset$.

\bigbreak

\ble\label{ch4:lem} Let $B, B^*\in\calb(\D)$ and define

\begin{align*}
\Omega( B,B^*) &:=\bigcap_{i=1}^d \big\{ \bigvee_{j\in\An(i): b_{j i}= b^*_{j i} } b_{j i} Z_j> \bigvee_{j\in\an(i):  b_{j i}\neq b^*_{j i}}(  b_{j i} \vee  b^                                                                                                                                                                                                                                                             *_{j i})  Z_j\big\}, \\  
A_{1/2}(B,B^*) & := \bigcap_{i=1}^d \big[  \big\{ \bfsx\in\R_+^d:x_i=\bigvee_{k\in\pa(i)}  b_{ki} x_k=\bigvee_{k\in\pa(i)} b^*_{ki}x_k \big \}\cup \big\{\bfsx\in\R_+^d: x_i>\bigvee_{k\in\pa(i)} ( b_{ki}\vee b^*_{ki})x_k \big \} \big].
\end{align*}
Then for every $F\in \calf$, 
\begin{align}\label{ch4:lem:eq3}
\P(F\cap \{  \bfx \in A_{1/2}(B,B^*) \})=\P(F\cap \Omega(B,B^*)).
\end{align}
\ele 
\bproof The proof is deferred to the appendix. 
\eproof

As a partition of $\R_+^d$ we now suggest $\big\{A_0(B,B^*), A_{1/2}(B,B^*), A_1(B,B^*)\big\}$, where $A_{1/2}(B, B^*)$ is defined above in Lemma~\ref{ch4:lem}, 
\[
A_0 (B,B^*) = \bigcup_{i\in V}\big[ \big\{ \bfsx\in\R_+^d:  x_i< \bigvee_{k\in\pa(i)}  b_{ki}x_k   \big\} \cup \big\{\bfsx\in\R_+^d: x_{i}  = \bigvee_{k\in \pa(i)} b^*_{ki}x_k >\bigvee_{k\in \pa(i)}  b_{k i}x_k  \big\}\big],\]
and $ A_1(B,B^*) =\R_+^d\setminus \big(A_0(B,B^*) \cup  A_{1/2}(B,B^*)\big)$.
With this partition we then have:
\bthe\label{ch4:dens:general1} Let $B, B^*\in\calb(\D)$. 
Then the function $\rho:\R_+^d\to \{0,1/2,1 \}$ 
\begin{align}\label{ch4:dens:gen1}
\bfsx\mapsto\rho(\bfsx, B, B^*)= 
\frac{1}{2}\cdot\bone_{A_{1/2}(B,B^*)}(\bfsx)+ \bone_{A_1(B,B^*)}(\bfsx) =
\begin{cases}
0,  & \text{if $\bfsx\in A_0(B, B^*)$}, \\
\frac{1}{2},  & \text{if $\bfsx\in A_{1/2}(B, B^*)$}, \\
1,  & \text{if $\bfsx\in A_1(B, B^*)$,}
\end{cases}
\end{align}
is a density of $P_B$ with respect to $P_B+P_{B^*}$. 
\ethe
\bproof{See the appendix}.
\eproof

We observe an interesting relation between  the density \eqref{ch4:dens:gen1} for $\D$ and  corresponding densities for subgraphs of $\D$.

\bexam[Local densities $\rho_i$\label{ch4:exam:margins}]\\
Consider the \DAG s 
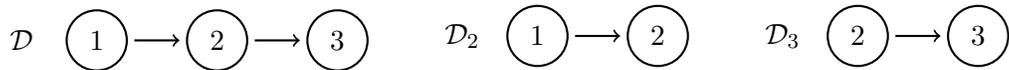
\begin{center}
	\setlength{\tabcolsep}{1em}
	\begin{tabular}{ccc}
		\begin{tikzpicture}[->,every node/.style={circle,draw},line width=0.8pt, node distance=1.6cm,minimum size=0.8cm,outer sep=1mm]
		\node (1) [outer sep=1mm]  {$1$};
		\node (2) [right of=1,outer sep=1mm] {$2$};
		\node (3) [right of=2,outer sep=1mm] {$3$}; 
		\foreach \from/\to in {2/3}
		\draw (\from) -- (\to);   
		\foreach \from/\to in {1/2}
		\draw (\from) -- (\to);  
		\node (n5)[draw=white,fill=white,left of=1,node distance=1cm] {$\D$}; 
		\end{tikzpicture}
		&
		\begin{tikzpicture}[->,every node/.style={circle,draw},line width=0.8pt, node distance=1.6cm,minimum size=0.8cm,outer sep=1mm]
		\node (1) [outer sep=1mm]  {$1$};
		\node (2) [right of=1,outer sep=1mm] {$2$}; 
		\foreach \from/\to in {1/2}
		\draw (\from) -- (\to);  
		\node (n5)[draw=white,fill=white,left of=1,node distance=1cm] {$\D_2$}; 
		\end{tikzpicture}
		&
		\begin{tikzpicture}[->,every node/.style={circle,draw},line width=0.8pt, node distance=1.6cm,minimum size=0.8cm,outer sep=1mm]
		\node (2) [right of=1,outer sep=1mm] {$2$};
		\node (3) [right of=2,outer sep=1mm] {$3$}; 
		\foreach \from/\to in {2/3}
		\draw (\from) -- (\to);    
		\node (n5)[draw=white,fill=white,left of=2,node distance=1cm] {$\D_3$}; 
		\end{tikzpicture}  
	\end{tabular}
\end{center}
Let  $\rho$, $\rho_2$, and $\rho_3$ be  the corresponding densities from \eqref{ch4:dens:gen1}. 
For the \ML\ coefficient matrix $B$ of a recursive \ML\ model on $\D$, let $B_{2}$ and $B_{3}$ be the \ML\ coefficient matrices of recursive \ML\ models on $\D_2$ and $\D_3$ with edge weight $c_{12}=b_{12}$ and $c_{23}=b_{23}$, and let starred quantities denote the same for $B^*$. 
We then find for $\bfsx=(x_1,x_2,x_3)\in\R_+^3$,
\begin{align*}
\lefteqn{ \rho(\bfsx,B,B^*) }\\
&=
\big(\rho_2(\bfsx_{\Pa(2)},B_{2},B^*_{2}) \vee \rho_3(\bfsx_{\Pa(3)},B_{3},B^*_{3})\big) \bone_{(0,\infty)}\big(   \rho_2(\bfsx_{\Pa(2)},B_{2},B^*_{2})\wedge  \rho_3(\bfsx_{\Pa(3)},B_{3},B^*_{3})\big).
\end{align*}
This can be observed from Figure~\ref{ch4:ex:fig:mar}, where 
the densities are depicted as  functions of  ${x_2}/{x_1}$ and/or ${x_3}/{x_2}$ for all nine different orders between the \ML\ coefficients in $B$ and $B^*$.
\begin{figure}[htb]
	\centering
	{\setlength{\tabcolsep}{0.0em}
		\renewcommand{\arraystretch}{0.5}
		\begin{tabular}{ccc}
			\begin{tikzpicture}[scale=0.55]
			\begin{scope}
			\clip[postaction={fill=red!40}](0,0)--(0,3.75)--(2.5,3.75)--(2.5,0);
			\end{scope}
			\begin{scope}
			\clip[postaction={fill=red!40}](2.5,0)--(2.5,2.5)--(3.75,3.75)--(3.75,0);
			\end{scope}
			\begin{scope}
			\clip[postaction={fill=OliveGreen!55}](2.5,1.25)--(2.5,2.5)--(3.75,2.5)--(3.75,1.25);
			\end{scope}
			\begin{scope}
			\clip[postaction={fill=blue!40}](2.5,2.5)--(2.5,3.75)--(3.75,3.75)--(3.75,2.5);
			\end{scope}
			\begin{scope}
			\clip[postaction={fill=OliveGreen!55}](1.25,1.25)--(1.25,2.5)--(2.5,2.5)--(2.5,1.25);
			\end{scope}
			\begin{scope}
			\clip[postaction={fill=OliveGreen!55}](1.25,2.5)--(1.25,3.75)--(2.5,3.75)--(2.5,2.5);
			\end{scope}
			\draw [line width=3pt,red!40] (2.5, 2.5) -- (2.5,3.75);
			\draw [line width=0.5pt,red!60] (0, 2.5) -- (2.5,2.5);
			\draw [line width=0.5pt,red!60] (0, 1.25) -- (2.5,1.25);
			\draw [line width=3pt,OliveGreen!55] (2.5, 1.25) -- (3.75,1.25);
			\draw [line width=0.5pt,red!60] (2.5,0) -- (2.5,1.25);
			\draw [line width=3pt,OliveGreen!55] (2.5,1.25) -- (2.5,2.5);
			\draw [line width=3pt,red!40] (2.5, 2.5) -- (3.75,2.5);
			\draw [line width=3pt,OliveGreen!55] (1.25, 2.5) -- (1.25,3.75);
			\draw [line width=3pt,red!40] (1.25, 2.5) -- (2.5,2.5);
			\draw [line width=3pt,red!40] (2.5, 1.25) -- (2.5,2.5);
			\draw [line width=3pt,OliveGreen!55] (1.25, 1.25) -- (2.5,1.25);
			\draw [line width=3pt,OliveGreen!55] (1.25, 1.25) -- (1.25,2.5);
			\draw [line width=0.5pt,red!60] (1.25,0) -- (1.25,1.25);
			\draw [red!40, fill=red!40]  (2.5,2.5) circle (4pt) node {};
			\draw [red!40, fill=red!40]  (2.5,1.25) circle (4pt) node {};
			\draw [red!40, fill=red!40]  (1.25,2.5) circle (4pt) node {};
			\draw [OliveGreen!55, fill=OliveGreen!55]  (1.25,1.25) circle (4pt) node {};
			\draw (2.5,-1.35) -- (2.5,-1.25) node[below=2pt] {\begin{scriptsize}
				$b^*_{12}$
				\end{scriptsize}};   
			\draw (-1.35,2.5) -- (-1.25,2.5) node[left=2pt] {\begin{scriptsize}
				$b^*_{23}$
				\end{scriptsize}};   
			\draw (1.25,-1.35) -- (1.25,-1.25) node[below=2pt] {\begin{scriptsize}
				$ b_{12}$
				\end{scriptsize}};   
			\draw (-1.35,1.25) -- (-1.25,1.25) node[left=2pt] {\begin{scriptsize}
				$ b_{23}$
				\end{scriptsize}};
			
			\coordinate (y) at (-1.25,4.5);
			\coordinate (x) at (-1.75,-1.25);
			\draw[<-] (y) node[left] {\begin{small}$\frac{x_3}{x_2}$\end{small}} -- (-1.25,-1.75);
			\draw[->] (x)-- (4.5,-1.25) node[below] {\begin{small}$\frac{x_2}{x_1}$\end{small}};        
			
			\draw [line width=3pt,red!40] (0,-0.625) -- (1.25,-0.625);
			\draw [line width=3pt,OliveGreen!55] (1.25,-0.625) -- (2.5,-0.625);
			\draw [line width=3pt,blue!40] (2.5,-0.625) -- (3.75,-0.625);
			
			\draw [OliveGreen!55, fill=OliveGreen!55]  (1.25,-0.625) circle (4pt) node {};           
			\draw [red!40, fill=red!40]  (2.5,-0.625) circle (4pt) node {};           
			\draw (3.75,-0.625) node [right] {\begin{scriptsize}$\rho_2$\end{scriptsize}};        
			
			\draw [line width=3pt,red!40] (-0.625,0) -- (-0.625,1.25);
			\draw [line width=3pt,OliveGreen!55] (-0.625,1.25) -- (-0.625,2.5);
			\draw [line width=3pt,blue!40] (-0.625,2.5) -- (-0.625,3.75);
			
			\draw [red!40, fill=red!40]  (-0.625,2.5) circle (4pt) node {};           
			\draw [OliveGreen!55, fill=OliveGreen!55]  (-0.625,1.25) circle (4pt) node {};       
			\draw (-0.625,3.75) node [above] {\begin{scriptsize}$\rho_3$\end{scriptsize}};   
			
			\draw (1.875,3.75) node [above=1pt] {\begin{scriptsize}$\rho(\cdot, B,B^*)$\end{scriptsize}}; 
			\end{tikzpicture}
			&

			\begin{tikzpicture}[scale=0.55]
			\begin{scope}
			\clip[postaction={fill=red!40}](0,0)--(0,3.75)--(1.875,3.75)--(1.875,0);
			\end{scope}
			\begin{scope}
			\clip[postaction={fill=red!40}](1.875,0)--(1.875,2.5)--(3.75,3.75)--(3.75,0);
			\end{scope}
			\begin{scope}
			\clip[postaction={fill=OliveGreen!55}](1.875,1.25)--(1.875,2.5)--(3.75,2.5)--(3.75,1.25);
			\end{scope}
			\begin{scope}
			\clip[postaction={fill=blue!40}](1.875,2.5)--(1.875,3.75)--(3.75,3.75)--(3.75,2.5);
			\end{scope}
			\draw [line width=3pt,OliveGreen!55] (1.875, 2.5) -- (1.875,3.75);
			\draw [line width=0.5pt,red!60] (0, 2.5) -- (1.875,2.5);
			\draw [line width=0.5pt,red!60] (0, 1.25) -- (1.875,1.25);
			\draw [line width=3pt,OliveGreen!55] (1.875, 1.25) -- (3.75,1.25);
			\draw [line width=0.5pt,red!60] (1.875,0) -- (1.875,1.25);
			\draw [line width=3pt,OliveGreen!55] (1.875,1.25) -- (1.875,2.5);
			\draw [line width=3pt,red!40] (1.875, 2.5) -- (3.75,2.5);
			\draw [line width=3pt,blue!40] (1.875, 2.5) -- (1.875,3.75);
			\draw [red!40, fill=red!40]  (1.875,2.5) circle (4pt) node {};
			\draw [OliveGreen!55, fill=OliveGreen!55]  (1.875,1.25) circle (4pt) node {};
			\draw (1.875,-1.35) -- (1.875,-1.25) node[below=2pt] {\begin{scriptsize}$ b_{12}=b^*_{12}$\end{scriptsize}};   
			\draw (-1.35,2.5) -- (-1.25,2.5) node[left=2pt] {\begin{scriptsize}$ b^*_{23}$\end{scriptsize}};      
			\draw (-1.35,1.25) -- (-1.25,1.25) node[left=2pt] {\begin{scriptsize}$ b_{23}$\end{scriptsize}};
			
			\coordinate (y) at (-1.25,4.5);
			\coordinate (x) at (-1.75,-1.25);
			\draw[<-] (y) node[left] {\begin{small}$\frac{x_3}{x_2}$\end{small}} -- (-1.25,-1.75);
			\draw[->] (x)-- (4.5,-1.25) node[below] {\begin{small}$\frac{x_2}{x_1}$\end{small}};        
			
			\draw [line width=3pt,red!40] (0,-0.625) -- (2.5,-0.625);
			\draw [line width=3pt,blue!40] (1.875,-0.625) -- (3.75,-0.625);

			\draw [blue!40, fill=blue!40]  (1.875,-0.625) circle (4pt) node {};           
			\draw (3.75,-0.625) node [right] {\begin{scriptsize}$\rho_2$\end{scriptsize}};        
			
			\draw [line width=3pt,red!40] (-0.625,0) -- (-0.625,1.25);
			\draw [line width=3pt,OliveGreen!55] (-0.625,1.25) -- (-0.625,2.5);
			\draw [line width=3pt,blue!40] (-0.625,2.5) -- (-0.625,3.75);
			
			\draw [red!40, fill=red!40]  (-0.625,2.5) circle (4pt) node {};           
			\draw [OliveGreen!55, fill=OliveGreen!55]  (-0.625,1.25) circle (4pt) node {};       
			
			\draw (-0.625,3.75) node [above] {\begin{scriptsize}$\rho_3$\end{scriptsize}};   
			
			\draw (1.875,3.75) node [above=1pt] {\begin{scriptsize}$\rho(\cdot, B,B^*)$\end{scriptsize}}; 
			
			\end{tikzpicture}

			&
			
			\begin{tikzpicture}[scale=0.55]
			\begin{scope}
			\clip[postaction={fill=red!40}](0,0)--(0,3.75)--(2.5,3.75)--(2.5,0);
			\end{scope}
			\begin{scope}
			\clip[postaction={fill=red!40}](2.5,0)--(2.5,2.5)--(3.75,3.75)--(3.75,0);
			\end{scope}
			\begin{scope}
			\clip[postaction={fill=OliveGreen!55}](2.5,1.25)--(2.5,2.5)--(3.75,2.5)--(3.75,1.25);
			\end{scope}
			\begin{scope}
			\clip[postaction={fill=blue!40}](2.5,2.5)--(2.5,3.75)--(3.75,3.75)--(3.75,2.5);
			\end{scope}
			\draw [line width=0.5pt,red!60] (1.25, 0) -- (1.25,3.75);
			\draw [line width=3pt,OliveGreen!55] (2.5, 2.5) -- (2.5,3.75);
			\draw [line width=0.5pt,red!60] (0, 2.5) -- (2.5,2.5);
			\draw [line width=0.5pt,red!60] (0, 1.25) -- (2.5,1.25);
			\draw [line width=3pt,OliveGreen!55] (2.5, 1.25) -- (3.75,1.25);
			\draw [line width=0.5pt,red!60] (2.5,0) -- (2.5,1.25);
			\draw [line width=3pt,OliveGreen!55] (2.5,1.25) -- (2.5,2.5);
			\draw [line width=3pt,red!40] (2.5, 2.5) -- (3.75,2.5);
			\draw [red!40, fill=red!40]  (2.5,2.5) circle (4pt) node {};
			\draw [OliveGreen!55, fill=OliveGreen!55]  (2.5,1.25) circle (4pt) node {};
			\draw (2.5,-1.35) -- (2.5,-1.25) node[below=2pt] {\begin{scriptsize}$ b_{12}$\end{scriptsize}};   
			\draw (-1.35,2.5) -- (-1.25,2.5) node[left=2pt] {\begin{scriptsize}$ b^*_{23}$\end{scriptsize}};   
			\draw (1.25,-1.35) -- (1.25,-1.25) node[below=2pt] {\begin{scriptsize}$ b^*_{12}$\end{scriptsize}};   
			\draw (-1.35,1.25) -- (-1.25,1.25) node[left=2pt] {\begin{scriptsize}$ b_{23}$\end{scriptsize}};
			\coordinate (y) at (-1.25,4.5);
			\coordinate (x) at (-1.75,-1.25);
			\draw[<-] (y) node[left] {\begin{small}$\frac{x_3}{x_2}$\end{small}} -- (-1.25,-1.75);
			\draw[->] (x)-- (4.5,-1.25) node[below] {\begin{small}$\frac{x_2}{x_1}$\end{small}};        
			
			\draw [line width=3pt,red!40] (0,-0.625) -- (2.5,-0.625);
			\draw [line width=3pt,blue!40] (2.5,-0.625) -- (3.75,-0.625);
			
			\draw [OliveGreen!60, fill=OliveGreen!60]  (2.5,-0.625) circle (4pt) node {};           
			\draw (3.75,-0.625) node [right] {\begin{scriptsize}$\rho_2$\end{scriptsize}};        
			
			\draw [line width=3pt,red!40] (-0.625,0) -- (-0.625,1.25);
			\draw [line width=3pt,OliveGreen!55] (-0.625,1.25) -- (-0.625,2.5);
			\draw [line width=3pt,blue!40] (-0.625,2.5) -- (-0.625,3.75);
			
			\draw [red!40, fill=red!40]  (-0.625,2.5) circle (4pt) node {};           
			\draw [OliveGreen!55, fill=OliveGreen!55]  (-0.625,1.25) circle (4pt) node {};                
			\draw (-0.625,3.75) node [above] {\begin{scriptsize}$\rho_3$\end{scriptsize}};   
			
			\draw (1.875,3.75) node [above=1pt] {\begin{scriptsize}$\rho(\cdot, B,B^*)$\end{scriptsize}};
			
			\end{tikzpicture}
			\\
			\begin{tikzpicture}[scale=0.55]   
			\begin{scope}
			\clip[postaction={fill=red!40}](0,0)--(0,3.75)--(2.5,3.75)--(2.5,0);
			\end{scope}
			\begin{scope}
			\clip[postaction={fill=red!40}](2.5,0)--(2.5,1.875)--(3.75,3.75)--(3.75,0);
			\end{scope}
			\begin{scope}
			\clip[postaction={fill=blue!40}](2.5,1.875)--(2.5,3.75)--(3.75,3.75)--(3.75,1.875);
			\end{scope}
			\begin{scope}
			\clip[postaction={fill=OliveGreen!55}](1.25,1.875)--(1.25,3.75)--(2.5,3.75)--(2.5,1.875);
			\end{scope}
			\draw [line width=3pt,red!40] (2.5,1.875) -- (2.5,3.75);
			\draw [line width=3pt,OliveGreen!55] (1.25, 1.875) -- (1.25,3.75);
			\draw [line width=0.5pt,red!60] (1.25, 0) -- (1.25,1.875);
			\draw [line width=0.5pt,red!60] (0, 1.875) -- (2.5,1.875);
			\draw [line width=0.5pt,red!60] (2.5,0) -- (2.5,1.875);
			\draw [line width=3pt,blue!40] (2.5, 1.875) -- (3.75,1.875);
			\draw [line width=3pt,OliveGreen!55] (1.25,1.875) -- (2.5,1.875);
			\draw[red!40, fill=red!40]  (2.5,1.875) circle (4pt) node {};
			\draw [OliveGreen!55, fill=OliveGreen!55]  (1.25,1.875) circle (4pt) node {};
			\draw (2.5,-1.35) -- (2.5,-1.25) node[below=2pt] {\begin{scriptsize}$ b^*_{12}$\end{scriptsize}};   
			\draw (-1.35,1.875) -- (-1.25,1.875) node [above=2pt,rotate=90pt] {\begin{scriptsize}$ b_{23}=b^*_{23}$\end{scriptsize}};   
			\draw (1.25,-1.35) -- (1.25,-1.25) node[below=2pt] {\begin{scriptsize}$  b_{12}$\end{scriptsize}};
			\coordinate (y) at (-1.25,4.5);
			\coordinate (x) at (-1.75,-1.25);
			\draw[<-] (y) node[left] {\begin{small}$\frac{x_3}{x_2}$\end{small}} -- (-1.25,-1.75);
			\draw[->] (x)-- (4.5,-1.25) node[below] {\begin{small}$\frac{x_2}{x_1}$\end{small}};        
			
			\draw [line width=3pt,red!40] (0,-0.625) -- (1.25,-0.625);
			\draw [line width=3pt,OliveGreen!55] (1.25,-0.625) -- (2.5,-0.625);
			\draw [line width=3pt,blue!40] (2.5,-0.625) -- (3.75,-0.625);
			
			\draw [OliveGreen!55, fill=OliveGreen!55]  (1.25,-0.625) circle (4pt) node {};           
			\draw [red!40, fill=red!40]  (2.5,-0.625) circle (4pt) node {};           
			\draw (3.75,-0.625) node [right] {\begin{scriptsize}$\rho_2$\end{scriptsize}};        
			
			\draw [line width=3pt,red!40] (-0.625,0) -- (-0.625,1.875);
			\draw [line width=3pt,blue!40] (-0.625,1.875) -- (-0.625,3.75);
			
			\draw [blue!40, fill=blue!40]  (-0.625,1.875) circle (4pt) node {};           
			\draw (-0.625,3.75) node [above] {\begin{scriptsize}$\rho_3$\end{scriptsize}};   
			
			\draw (1.875,3.75) node [above=1pt] {\begin{scriptsize}$\rho(\cdot, B,B^*)$\end{scriptsize}}; 
			\end{tikzpicture}			
			&
			\begin{tikzpicture}[scale=0.55]  
			\begin{scope}
			\clip[postaction={fill=red!40}](0,0)--(0,3.75)--(1.875,3.75)--(1.875,0);
			\end{scope}
			\begin{scope}
			\clip[postaction={fill=red!40}](1.875,0)--(1.875,1.875)--(3.75,3.75)--(3.75,0);
			\end{scope}
			\begin{scope}
			\clip[postaction={fill=blue!40}](1.875,1.875)--(1.875,3.75)--(3.75,3.75)--(3.75,1.875);
			\end{scope}
			\draw [line width=3pt,blue!40] (1.875, 1.875) -- (1.875,3.75);
			\draw [line width=0.5pt,red!60] (0, 1.875) -- (1.875,1.875);
			\draw [line width=0.5pt,red!60] (1.875,0) -- (1.875,1.875);
			\draw [line width=3pt,blue!40] (1.875,1.875) -- (3.75,1.875);
			\draw [blue!40, fill=blue!40]  (1.875,1.875) circle (4pt) node {};
			\draw (1.875,-1.35) -- (1.875,-1.25) node[below=2pt] {\begin{scriptsize}$ b_{12}=b^*_{12}$\end{scriptsize}};   
			\draw (-1.35,1.875) -- (-1.25,1.875) node [above=2pt,rotate=90pt] {\begin{scriptsize}$ b_{23}=b^*_{23}$\end{scriptsize}}; 
			\coordinate (y) at (-1.25,4.5);
			\coordinate (x) at (-1.75,-1.25);
			\draw[<-] (y) node[left] {\begin{small}$\frac{x_3}{x_2}$\end{small}} -- (-1.25,-1.75);
			\draw[->] (x)-- (4.5,-1.25) node[below] {\begin{small}$\frac{x_2}{x_1}$\end{small}};        
			
			\draw [line width=3pt,red!40] (0,-0.625) -- (2.5,-0.625);
			\draw [line width=3pt,blue!40] (1.875,-0.625) -- (3.75,-0.625);

			\draw [blue!40, fill=blue!40]  (1.875,-0.625) circle (4pt) node {};           
			\draw (3.75,-0.625) node [right] {\begin{scriptsize}$\rho_2$\end{scriptsize}};        
			
			\draw [line width=3pt,red!40] (-0.625,0) -- (-0.625,1.875);
			\draw [line width=3pt,blue!40] (-0.625,1.875) -- (-0.625,3.75);
			
			\draw [blue!40, fill=blue!40]  (-0.625,1.875) circle (4pt) node {};           
			\draw (-0.625,3.75) node [above] {\begin{scriptsize}$\rho_3$\end{scriptsize}};
			
			\draw (1.875,3.75) node [above=1pt] {\begin{scriptsize}$\rho(\cdot, B,B^*)$\end{scriptsize}}; 
			\end{tikzpicture}
			&		
			\begin{tikzpicture}[scale=0.55] 
			\begin{scope}
			\clip[postaction={fill=red!40}](0,0)--(0,3.75)--(2.5,3.75)--(2.5,0);
			\end{scope}
			\begin{scope}
			\clip[postaction={fill=red!40}](2.5,0)--(2.5,1.875)--(3.75,3.75)--(3.75,0);
			\end{scope}
			\begin{scope}
			\clip[postaction={fill=blue!40}](2.5,1.875)--(2.5,3.75)--(3.75,3.75)--(3.75,1.875);
			\end{scope}
			\draw [line width=0.5pt,red!60] (1.25, 0) -- (1.25,3.75);
			\draw [line width=3pt,OliveGreen!55] (2.5, 1.875) -- (2.5,3.75);
			\draw [line width=0.5pt,red!60] (0, 1.875) -- (2.5,1.875);
			\draw [line width=0.5pt,red!60] (2.5,0) -- (2.5,1.875);
			\draw [line width=3pt,blue!40] (2.5, 1.875) -- (3.75,1.875);
			
			\draw [OliveGreen!55, fill=OliveGreen!55]  (2.5,1.875) circle (4pt) node {};
			\draw (2.5,-1.35) -- (2.5,-1.25) node[below=2pt] {\begin{scriptsize}$ b_{12}$\end{scriptsize}};   
			\draw (-1.35,1.875) -- (-1.25,1.875) node [above=2pt,rotate=90pt] {\begin{scriptsize}$ b_{23}=b^*_{23}$\end{scriptsize}};   
			\draw (1.25,-1.35) -- (1.25,-1.25) node[below=2pt] {\begin{scriptsize}$ b^*_{12}$\end{scriptsize}}; 
			
			\coordinate (y) at (-1.25,4.5);
			\coordinate (x) at (-1.75,-1.25);
			\draw[<-] (y) node[left] {\begin{small}$\frac{x_3}{x_2}$\end{small}} -- (-1.25,-1.75);
			\draw[->] (x)-- (4.5,-1.25) node[below] {\begin{small}$\frac{x_2}{x_1}$\end{small}};        
			
			\draw [line width=3pt,red!40] (0,-0.625) -- (2.5,-0.625);
			\draw [line width=3pt,blue!40] (2.5,-0.625) -- (3.75,-0.625);
			
			\draw [OliveGreen!60, fill=OliveGreen!60]  (2.5,-0.625) circle (4pt) node {};           
			\draw (3.75,-0.625) node [right] {\begin{scriptsize}$\rho_2$\end{scriptsize}};        
			
			\draw [line width=3pt,red!40] (-0.625,0) -- (-0.625,1.875);
			\draw [line width=3pt,blue!40] (-0.625,1.875) -- (-0.625,3.75);
			
			\draw [blue!40, fill=blue!40]  (-0.625,1.875) circle (4pt) node {};           
			\draw (-0.625,3.75) node [above] {\begin{scriptsize}$\rho_3$\end{scriptsize}};   
			
			\draw (1.875,3.75) node [above=1pt] {\begin{scriptsize}$\rho(\cdot, B,B^*)$\end{scriptsize}};
			
			\end{tikzpicture}

			\\
			
			\begin{tikzpicture}[scale=0.55]
			\begin{scope}
			\clip[postaction={fill=red!40}](0,0)--(0,3.75)--(2.5,3.75)--(2.5,0);
			\end{scope}
			\begin{scope}
			\clip[postaction={fill=red!40}](2.5,0)--(2.5,2.5)--(3.75,3.75)--(3.75,0);
			\end{scope}
			\begin{scope}
			\clip[postaction={fill=blue!40}](2.5,2.5)--(2.5,3.75)--(3.75,3.75)--(3.75,2.5);
			\end{scope}
			\begin{scope}
			\clip[postaction={fill=OliveGreen!55}](1.25,2.5)--(1.25,3.75)--(2.5,3.75)--(2.5,2.5);
			\end{scope}
			
			\draw [line width=0.5pt,red!60] (0, 1.25) -- (3.75,1.25);
			\draw [line width=3pt,red!40] (2.5, 2.5) -- (2.5,3.75);
			\draw [line width=0.5pt,red!60] (0, 2.5) -- (2.5,2.5);
			\draw [line width=0.5pt,red!60] (2.5,0) -- (2.5,2.5);
			\draw [line width=3pt,OliveGreen!55] (2.5, 2.5) -- (3.75,2.5);
			\draw [line width=3pt,OliveGreen!55] (1.25, 2.5) -- (1.25,3.75);
			\draw [line width=0.5pt,red!60] (1.25, 0) -- (1.25,2.5);
			\draw [line width=3pt,OliveGreen!55] (1.25, 2.5) -- (2.5,2.5);

			\draw [red!40, fill=red!40]  (2.5,2.5) circle (4pt) node {};
			\draw [OliveGreen!55, fill=OliveGreen!55]  (1.25,2.5) circle (4pt) node {};
			\draw (2.5,-1.35) -- (2.5,-1.15) node[below=2pt] {\begin{scriptsize}$ b^*_{12}$\end{scriptsize}};   
			\draw (-1.35,2.5) -- (-1.15,2.5) node[left=2pt] {\begin{scriptsize}$ b_{23}$\end{scriptsize}};   
			\draw (1.25,-1.35) -- (1.25,-1.15) node[below=2pt] {\begin{scriptsize}$  b_{12}$\end{scriptsize}};   
			\draw (-1.35,1.25) -- (-1.15,1.25) node[left=2pt] {\begin{scriptsize}
				$ b^*_{23}$
				\end{scriptsize}};
			
			\coordinate (y) at (-1.25,4.5);
			\coordinate (x) at (-1.75,-1.25);
			\draw[<-] (y) node[left] {\begin{small}$\frac{x_3}{x_2}$\end{small}} -- (-1.25,-1.75);
			\draw[->] (x)-- (4.5,-1.25) node[below] {\begin{small}$\frac{x_2}{x_1}$\end{small}};        
			
			\draw [line width=3pt,red!40] (0,-0.625) -- (1.25,-0.625);
			\draw [line width=3pt,OliveGreen!55] (1.25,-0.625) -- (2.5,-0.625);
			\draw [line width=3pt,blue!40] (2.5,-0.625) -- (3.75,-0.625);
			
			\draw [OliveGreen!55, fill=OliveGreen!55]  (1.25,-0.625) circle (4pt) node {};           
			\draw [red!40, fill=red!40]  (2.5,-0.625) circle (4pt) node {};           
			\draw (3.75,-0.625) node [right] {\begin{scriptsize}$\rho_2$\end{scriptsize}};        
			
			\draw [line width=3pt,red!40] (-0.625,0) -- (-0.625,2.5);
			\draw [line width=3pt,blue!40] (-0.625,2.5) -- (-0.625,3.75);
			
			\draw [OliveGreen!60, fill=OliveGreen!60]  (-0.625,2.5) circle (4pt) node {};           
			\draw (-0.625,3.75) node [above] {\begin{scriptsize}$\rho_3$\end{scriptsize}};   
			
			\draw (1.875,3.75) node [above=1pt] {\begin{scriptsize}$\rho(\cdot, B,B^*)$\end{scriptsize}}; 
			
			\end{tikzpicture}

			&
			
			\begin{tikzpicture}[scale=0.55]
			\begin{scope}
			\clip[postaction={fill=red!40}](0,0)--(0,3.75)--(1.875,3.75)--(2.5,0);
			\end{scope}
			\begin{scope}
			\clip[postaction={fill=red!40}](1.875,0)--(1.875,2.5)--(3.75,3.75)--(3.75,0);
			\end{scope}
			\begin{scope}
			\clip[postaction={fill=blue!40}](1.875,2.5)--(1.875,3.75)--(3.75,3.75)--(3.75,2.5);
			\end{scope}
			
			\draw [line width=0.5pt,red!60] (0,1.25) -- (3.75,1.25);
			\draw [line width=3pt,blue!40] (1.875, 2.5) -- (1.875,3.75);
			\draw [line width=0.5pt,red!60] (0, 2.5) -- (1.875,2.5);
			\draw [line width=0.5pt,red!60] (1.875,0) -- (1.875,2.5);
			\draw [line width=3pt,OliveGreen!55] (1.875, 2.5) -- (3.75,2.5);
			\draw [OliveGreen!55, fill=OliveGreen!55]  (1.875,2.5) circle (4pt) node {};
			\draw (1.875,-1.35) -- (1.875,-1.15) node[below=2pt] {\begin{scriptsize}$ b_{12}=b^*_{12}$\end{scriptsize}};   
			\draw (-1.35,2.5) -- (-1.15,2.5) node[left=2pt] {\begin{scriptsize}$ b_{23}$\end{scriptsize}};    
			\draw (-1.35,1.25) -- (-1.25,1.25) node[left=2pt] {\begin{scriptsize}$ b^*_{23}$\end{scriptsize}};

			\coordinate (y) at (-1.25,4.5);
			\coordinate (x) at (-1.75,-1.25);
			\draw[<-] (y) node[left] {\begin{small}$\frac{x_3}{x_2}$\end{small}} -- (-1.25,-1.75);
			\draw[->] (x)-- (4.5,-1.25) node[below] {\begin{small}$\frac{x_2}{x_1}$\end{small}};        
			
			\draw [line width=3pt,red!40] (0,-0.625) -- (2.5,-0.625);
			\draw [line width=3pt,blue!40] (1.875,-0.625) -- (3.75,-0.625);

			\draw [blue!40, fill=blue!40]  (1.875,-0.625) circle (4pt) node {};           
			\draw (3.75,-0.625) node [right] {\begin{scriptsize}$\rho_2$\end{scriptsize}};        
			
			\draw [line width=3pt,red!40] (-0.625,0) -- (-0.625,2.5);
			\draw [line width=3pt,blue!40] (-0.625,2.5) -- (-0.625,3.75);
			
			\draw [OliveGreen!60, fill=OliveGreen!60]  (-0.625,2.5) circle (4pt) node {};           
			\draw (-0.625,3.75) node [above] {\begin{scriptsize}$\rho_3$\end{scriptsize}};   
			
			\draw (1.875,3.75) node [above=1pt] {\begin{scriptsize}$\rho(\cdot, B,B^*)$\end{scriptsize}}; 
			
			\end{tikzpicture}

			&
			\begin{tikzpicture}[scale=0.55]
			\begin{scope}
			\clip[postaction={fill=red!40}](0,0)--(0,3.75)--(2.5,3.75)--(2.5,0);
			\end{scope}
			\begin{scope}
			\clip[postaction={fill=red!40}](2.5,0)--(2.5,2.5)--(3.75,3.75)--(3.75,0);
			\end{scope}
			\begin{scope}
			\clip[postaction={fill=blue!40}](2.5,2.5)--(2.5,3.75)--(3.75,3.75)--(3.75,2.5);
			\end{scope}
			\draw [line width=0.5pt,red!60] (1.25, 0) -- (1.25,3.75);
			\draw [line width=0.5pt,red!60] (0, 1.25) -- (3.75,1.25);
			\draw [line width=3pt,OliveGreen!55] (2.5, 2.5) -- (2.5,3.75);
			\draw [line width=0.5pt,red!60] (0, 2.5) -- (2.5,2.5);
			\draw [line width=0.5pt,red!60] (2.5,0) -- (2.5,2.5);
			\draw [line width=3pt,OliveGreen!55] (2.5, 2.5) -- (3.75,2.5);
			
			\draw [OliveGreen!55, fill=OliveGreen!55]  (2.5,2.5) circle (4pt) node {};
			\draw (2.5,-1.25) -- (2.5,-1.35) node[below=2pt] {\begin{scriptsize}{$ b_{12}$}\end{scriptsize}};   
			\draw (-1.35,2.5) -- (-1.25,2.5) node[left=2pt] {\begin{scriptsize}$b_{23}$\end{scriptsize}};   
			\draw (1.25,-1.35) -- (1.25,-1.25) node[below=2pt] {\begin{scriptsize}$ b^*_{12}$\end{scriptsize}};   
			\draw (-1.35,1.25) -- (-1.25,1.25) node[left=2pt] {\begin{scriptsize}$ b^*_{23}$\end{scriptsize}};
			
			\coordinate (y) at (-1.25,4.5);
			\coordinate (x) at (-1.75,-1.25);
			\draw[<-] (y) node[left] {\begin{small}$\frac{x_3}{x_2}$\end{small}} -- (-1.25,-1.75);
			\draw[->] (x)-- (4.5,-1.25) node[below] {\begin{small}$\frac{x_2}{x_1}$\end{small}};        
			
			\draw [line width=3pt,red!40] (0,-0.625) -- (2.5,-0.625);
			\draw [line width=3pt,blue!40] (2.5,-0.625) -- (3.75,-0.625);
			
			\draw [OliveGreen!60, fill=OliveGreen!60]  (2.5,-0.625) circle (4pt) node {};           
			\draw (3.75,-0.625) node [right] {\begin{scriptsize}$\rho_2$\end{scriptsize}};        
			
			\draw [line width=3pt,red!40] (-0.625,0) -- (-0.625,2.5);
			\draw [line width=3pt,blue!40] (-0.625,2.5) -- (-0.625,3.75);
			
			\draw [OliveGreen!60, fill=OliveGreen!60]  (-0.625,2.5) circle (4pt) node {};           
			\draw (-0.625,3.75) node [above] {\begin{scriptsize}$\rho_3$\end{scriptsize}};   
			
			\draw (1.875,3.75) node [above=1pt] {\begin{scriptsize}$\rho(\cdot, B,B^*)$\end{scriptsize}};

			\end{tikzpicture}
	\end{tabular}}
	\caption{The densities 
		$\rho(\bfsx, 
		B,B^*)$, $\rho_2(\bfsx_{\Pa(2)},
		B_{2},B^*_{2})$, $\rho_3(\bfsx_{\Pa(3)},
		B_{3},B^*_{3})$  
		from  Example~\ref{ch4:exam:margins} as functions of ${x_2}/{x_1}$ and/or ${x_3}/{x_2}$. The area where the density is ${\color{red!40}\boldsymbol{0}}/{\color{blue!40}\boldsymbol{\frac{1}{2}}}/{\color{OliveGreen!55}\boldsymbol{1}}$ is coloured in {\color{red!40}\bf red}/{\color{blue!40} \bf blue}/{\color{OliveGreen!55} \bf green}. 
		\label{ch4:ex:fig:mar} 	
	}
\end{figure} 

 Conversely,  $\rho_2$ and $\rho_3$ can be derived from $\rho$ as follows:
\begin{align*}
\rho_2(\bfsx_{\Pa(2)},B_{12},B^*_{12})&= \min_{ \{y\in\R_+: \rho((\bfsx_{\Pa(2)},y),B,B^*)>0\}} \rho((\bfsx_{\Pa(2)},y),B,B^*), \\ 
\rho_3(\bfsx_{\Pa(3)},B_{23},B^*_{23})&= \min_{\{y\in\R_+: \rho((y,\bfsx_{\Pa(3)}),B,B^*)>0\}} \rho((y,\bfsx_{\Pa(3)}),B,B^*),
\end{align*}which we learn from Figure~\ref{ch4:ex:fig:mar} again. 
\halmos
\eexam

We now extend the findings from  Example~\ref{ch4:exam:margins} to the general case. Furthermore, we show that the densities $\rho_i$ are densities of regular conditional distributions.  
\bpr\label{ch4:prop:rhoi}
Let $B, B^*\in\calb(\D)$ and let $\bfx= \bfz\odot B, \bfx^*=\bfz\odot B^*$ follow corresponding recursive \ML\ models on $\D$. 
For $i\in V$, let $\rho_i$ be the density given in \eqref{ch4:dens:gen1} with respect to the \DAG\ $\D_i=(\Pa(i),\{(k,i): k\in\pa(i)\})$ as well as $B_i$ and $B_i^*$ the \ML\ coefficient matrices of recursive \ML\ models on $\D_i$  with edge weights $c_{ki}=b_{ki}$ and $c^*_{ki}=b^*_{ki}$, respectively. 
\begin{enumerate}
	\item[(a)]
	We have for 
	$\rho(\bfsx,B,B^*)$ given in  \eqref{ch4:dens:gen1} 
		\begin{align}\label{ch4:rhorhoi}
	\rho(\bfsx,B,B^*)&= \big( \bigvee_{i\in V} \rho_i(\bfsx_{\Pa(i)},B_i,B^*_i)\big)  \bone_{(0,\infty)}\big(\bigwedge_{i\in V} \rho_i(\bfsx_{\Pa(i)},B_i,B^*_i)\big).
	\end{align}
	
	\item[(b)] The function $\rho_i$ can be computed from $\rho$ by
	\begin{align*}
	\rho_i(\bfsx_{\Pa(i)},B_i,B^*_i)= \min_{\{\bfsy\in \R_+^d:\bfsy_{\Pa(i)}=\bfsx_{\Pa(i)}, \rho(\bfsy,B,B^*)>0\}} \rho(\bfsy,B,B^*),
	\end{align*}
	where we set $\min_{\bfsy\in\emptyset}\rho(\bfsy,B,B^*)=0$. 
	\item[(c)] The function $\rho_i:\R_+^d\to\{0,1/2,1 \}$ such that  $\bfsx_{\Pa(i)}\mapsto\rho_i(\bfsx_{\Pa(i)},B_i,B_i^*)$ is  a density of $P^{i\mid\pa(i)}_B$ with respect to  $P^{i\mid\pa(i)}_B+P^{i\mid\pa(i)}_{B^*}$, where $P_B^{i\mid\pa(i)}$ is a regular conditional distribution of $X_i$ given $\bfx_{\pa(i)}$ and $P_{B^*}^{i\mid\pa(i)}$ one of $X^*_i$ given $\bfx^*_{\pa(i)}$. 
\end{enumerate}
\epr 
\bproof See the appendix. 
\eproof

Next, we show that $\widehat B$ is indeed a GMLE in the sense of \cite{kiefer1956}. Note also that the GMLE is obtained by piecing together individual GMLEs corresponding to conditional distributions of any variable given its parents.  Thus this is similar to what is obtained in cases where the distributions have densities with respect to a product measure, as the maximum of the likelihood function is then obtained by maximizing each conditional likelihood function for the density  of a node given its parents.
\bthe\label{ch4:the:gmle}  Let   $\bfsx^{(t)}=\big(x_1^{(t)},\ldots,x_n^{(t)}\big)$ for $t=1,\ldots,n$ be a sample from a recursive \ML\ model  on a \DAG\ $\D$ with \ML\ coefficient matrix
 $B\in \calb(\D)$ unknown. 
\begin{enumerate}
	\item[(a)] The matrix $\wh B$ from \eqref{ch4:eq:Bhat} is a GMLE of $B$.
	\item[(b)] For every $i\in V$,  
	$(\wh b_{ki},k\in\pa(i))$ 
	 is a GMLE of the \ML\ coefficients  $(b_{ki},k\in\pa(i))$  of a random vector following a  recursive ML model  on $\D_i=(\Pa(i),\{(k,i): k\in\pa(i)\})$ with edge weights $c_{ki}=b_{ki}$. 
	\item[(c)] For every $i\in V$ and $k\in\pa(i)$, $\wh b_{ki}$ is the only GMLE of the \ML\ coefficient $b_{ki}$ of a  random vector following a recursive \ML\ model  on  $\D_{ki}=(\{k,i\},\{(k,i)\})$ with edge weight  $c_{ki}=b_{ki}$. 
\end{enumerate} 
\ethe

\bproof
(a) First, recall that $\widehat B$ is indeed a \CM\ of a \mSEM\ on $\D$.  The first condition in the definition of a GMLE in \eqref{ch4:def:gmle} is satisfied due to the definition of  $\rho(\cdot,\wh B,\wh B)$ since $A_{1/2}(\wh B,\wh B)=\mathbb{R}_+^d$.  Since  the densities $\rho(\cdot,\wh B, B)$ and  $\rho(\cdot, B, \wh B)$ have the values $0$, $1$, $1/2$, and  $A_{1/2}(\wh B,B)=A_{1/2}(B,\wh B)$, to verify the second condition in \eqref{ch4:def:gmle}, it suffices to show that  there is some realization $\bfsx^{(t_1)}\in A_0(B,\wh B)$ whenever there is  some realization $\bfsx^{(t_2)}\in A_0(\widehat B, B)$; cf.~Example~\ref{ch4:ex1}, in particular (\ref{ch4:cond2:DAG12}). So  let $\bfsx^{(t_2)}\in A_0(\widehat B, B)$ for some $t_2\in\{1,\ldots, n\}$. 
We find, for some $i\in V$, from  the definition of $A_0(\widehat B, B)$ and  the fact that $x_i^{(t)} \ge \bigvee_{k\in\pa(i)} \wh b_{ki} x_k^{(t)}$,
\begin{align*}
\bfsx^{(t_2)}\in \big\{ \bfsx\in\R_+^d: \bigvee_{k\in\pa(i)} \widehat b_{ki}x_k < x_i=\bigvee_{k\in\pa(i)}  b_{ki}x_k\big\}.
\end{align*}
Hence, $x^{(t_2)}_i= b_{ki}x^{(t_2)}_k$ for some $k\in \pa(i)$ with $\widehat b_{ki}<  b_{ki}$. Let now $t_1\in\{1,\ldots,n\}$ such that $\bigwedge_{s=1}^n  y^{(s)}_{ki}=y^{(t_1)}_{ki}$. As $\widehat b_{ki}=\bigwedge_{s=1}^n y^{(s)}_{ki}$, we have  $x_i^{(t_1)}< b_{ki} x_k^{(t_1)}$ implying that $\bfsx^{(t_1)} \in A_0(B,\wh B)$. \\
The statement in (b) is a consequence of (a), and
(c) has already been shown in Example~\ref{ch4:ex1}.
\eproof

\begin{figure}[htb]
	\centering
	\begin{tikzpicture}[->,every node/.style={circle,draw},line width=0.8pt, node distance=1.6cm,minimum size=0.8cm,outer sep=1mm]
	\node (1)  {$1$};
	\node (2) [right of=1] {$2$}; 
	\node (3) [below of=2] {$3$};
	\node (4) [below of=1] {$4$};
	\node (n5)[draw=white,fill=white,left of=4,node distance=1cm] {$\D$};
	\foreach \from/\to in {1/2,1/4,2/3,2/4,3/4}
	\draw (\from) -- (\to);   
	\node[draw=none,below right=-0.4cm and 0.1cm of 4] {\footnotesize $b_{34}$}; 
	\node[draw=none,above right=0.1cm and 0.25cm of 4,rotate=45] 
	{\footnotesize $b_{24}$};
	\node[draw=none,below right=0cm and -0.4cm of 2] {\footnotesize $b_{23}$};  
	\node[draw=none,below left=0cm and -0.45cm of 1] {\footnotesize $b_{14}$};         
	\node[draw=none,above right=-0.45cm and 0cm of 1] {\footnotesize $b_{12}$};
	\node (11) [node distance=2cm,right of=2]  {$1$};
	\node (n5)[draw=white,fill=white,below of=11,node distance=1cm] {$\D_1$};   
	\node (12) [right of=11,node distance=2cm] {$1$};
	\node (22) [right of=12] {$2$}; 
	\node (42) [below of=12,color=white,node distance=1cm] {\color{black}$\D_2$};
	\node[draw=none,above right=-0.45cm and 0cm of 12] {\footnotesize  $b_{12}$};

	\foreach \from/\to in {12/22}
	\draw          (\from) -- (\to);   
	
	\node (23) [right of=22,node distance=2cm] {$2$};
	\node (33) [below of=23] {$3$};
	\node (13) [color=white,left of=33,node distance=1cm] {\color{black}$\D_3$};
	\foreach \from/\to in {23/33}
	\draw (\from) -- (\to);   
	\node[draw=none,below right=0cm and -0.4cm of 23] {\footnotesize $b_{23}$};   
	\node (14) [right of =23,node distance=2cm]  {$1$};
	\node (24) [right of=14] {$2$}; 
	\node (34) [below of=24] {$3$};
	\node (44) [below of=14] {$4$};
	\node (n5)[draw=white,fill=white,right of=34,node distance=1cm] {$\D_4$};
	\foreach \from/\to in {14/44,24/44,34/44}
	\draw (\from) -- (\to); 
	\node[draw=none,below left=0cm and -0.45cm of 14] {\footnotesize $b_{14}$}; 
	\node[draw=none,below right=-0.4cm and 0.1cm of 44] {\footnotesize $b_{34}$};   
	\node[draw=none,above right=0.1cm and 0.25cm of 44,rotate=45]  
	{\footnotesize $b_{24}$};
	\end{tikzpicture}
	\caption{The \DAG s $\D_i$ 
		 from Theorem~\ref{ch4:the:gmle}(b) 
	 for a recursive \ML\ model on  the \DAG\ $\D$  depicted on the left-hand side with \ML\  coefficient matrix $B$. The edges are marked with the corresponding \ML\ coefficients. Note that  $b_{12}, b_{14}, b_{34}, b_{24}$ can be arbitary positive numbers but $b_{24}\ge b_{23}b_{34}$.\label{ch4:exam:Di}}
\end{figure}
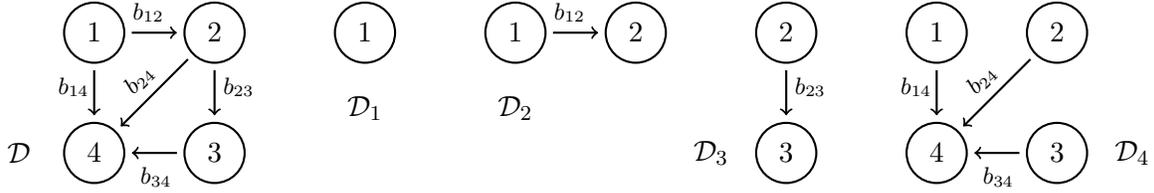

Figure~\ref{ch4:exam:Di} illustrates the \DAG s $\D_i$ in Theorem~\ref{ch4:the:gmle}(b) or Proposition~\ref{ch4:prop:rhoi}. 

\section{Learning the structure of a recursive max-linear model\label{ch4:s6}} 

In contrast to the assumptions in the previous section, we now assume  independent realizations  $\bfsx^{(1)},\ldots, \bfsx^{(n)}$ of $\bfx$ following a recursive \ML\ model but the underlying \DAG\ $\D$ is unknown. We know from  previous discussions that it is not possible to recover $\D$ and the true edge weights $c_{ki}$, and we therefore again focus on the estimation of  $B$.  

Following Algorithm~\ref{ch4:alg1}, it suffices  for any pair of distinct $i,j\in V$ to decide whether  $\ubpp(Y_{ji})= \ubpp(X_i/X_j)$  has a positive lower bound, alternatively a finite upper bound, and if so, to estimate the bound. Recall from Table~\ref{ch4:table:rel} that, if there is such a bound, then it is an atom of $Y_{ji}$. Since we can  expect to observe atoms more than twice for $n$ sufficiently large, we propose the following estimation method.

\balg[Find an estimate  $\widecheck B$ of $B$ from $\bfsx^{(1)},\ldots, \bfsx^{(n)}$\label{ch4:alg4}]
\begin{enumerate}
	\item[1.] For all $i\in V=\{1,\ldots, d\}$, set $\widecheck b_{ii}=1$. 
	\item[2.] For all $i,j\in V$ with $i\neq j$ ,\\[-0.95em]
	\begin{enumerate}
		\item[] if $\#\left\{t:\bigwedge_{s=1}^ny_{ji}^{(s)}=y_{ji}^{(t)}\right\}\ge 2$, then conclude $j\in \an(i)$, set $\widecheck{b}_{ji}=\bigwedge_{t=1}^n y_{ji}^{(t)}$;
		\begin{enumerate}
			\item[] \hphantom{ else,\,\,\,} else, set $\widecheck{b}_{ji}=0$.
		\end{enumerate}
	\end{enumerate}
\end{enumerate}
\ealg

The second item summarizes two steps: the first is concerned with estimating the ancestors of the nodes, 
 the second with estimating the \ML\ coefficients.
 
 Note that the estimate $\widecheck B$ from Algorithm~\ref{ch4:alg4} is not necessarily a ML coefficient matrix of a recursive ML model. For example, the property  that $b_{ji}>0$ if $b_{jk}b_{ki}>0$ (see, for example, Corollary~3.12 of \cite{GK1}) is not guaranteed. 
Many modifications of $\widecheck B$ are possible, and here we shall not discuss this in detail. Rather we notice that the probability that
 Algorithm~\ref{ch4:alg4} outputs the true  \ML\ coefficient matrix $B$ tends to one as $n\to \infty$. 
  As in the case where the DAG is known --- see  Proposition~\ref{ch4:PropBhat} ---  this probability 
  converges to one at an exponential rate.

\section{Conclusion and outlook\label{ch4:s7}} 
 
We studied the identifiability of the elements of a recursive \ML\ model from the distribution $\call(\bfx)$ of $\bfx$. The associated \DAG\ and the   edge weights are not identifiable, however, the \ML\ coefficient matrix $B$ is. In other words, we can identify the representation \eqref{ch4:ml-noise} but not \eqref{ch4:ml-sem}. 
The class of all \DAG s and edge weights  that could have generated $\bfx$ via  \eqref{ch4:ml-sem} and the distribution of the innovation vector 
are  identifiable from $\call(\bfx)$.  As a consequence, we can recover $B$, the class of the DAGs and edge weights, and the innovation distributions  from realizations of $\bfx$. 

We have shown that $\hat B$ is a generalized maximum likelihood estimate. This is primarily of theoretical interest as it shows the estimate is not purely based on an \emph{ad hoc procedure}. However, it opens up the possibility of going further, using likelihood theory, for example to study issues of likelihood ratio testing of hypothesis for specific values of the coefficients, or even for the presence or absence of edges in the underlying graph. 


Parameter estimation and structure learning for recursive \ML\ models seem to be  challenging tasks because  assumptions usually made in standard methods are not met. However, in both cases,  $B$ can be estimated by a simple procedure.  The key idea of our approach is to consider the observed ratios between any pair of components, i.e.\  to perform a transformation on the realizations. The  transformed realizations or rather the distributional properties of the corresponding random variables  make it possible to identify, with probability 1, the true $B$  whenever  the number of observations $n$ is sufficiently large. It would be interesting to investigate the relationship between the performance of our procedures and the number $n$ of observations. Here, one possible question is how many observations are at least necessary to estimate $B$ exactly; see,  Example~\ref{ch4:examBhat}.  In addition it would be interesting to study estimation of the DAG structure for moderate sample sizes, where exact estimation is not guaranteed.

We emphasize again that, although our estimates are derived under the assumption that the distribution of the innovation vector $\bfz$  is fixed, the estimates do not depend on what this distribution is and would therefore also be valid in the situation where  the innovations are independent with unkown distributions that are atom-free and have  support equal to $\R_+$. 
Algorithm~\ref{ch4:alg3} provides a recursive procedure to obtain the distribution functions $F_{Z_i}$ from $B$ and  the marginal distribution functions $F_{X_i}$ of $X_i$.
Estimating $B$ by $\wh B$ and the distributions   $F_{X_i}$, for example, by their empirical versions, we can apply this procedure  to find estimators 
 of the distributions $F_{Z_i}$ although it will formally 
 violate  the assumption of atom-freeness and thus it is both more efficient and formally correct to estimate these parametrically, or under suitable monotonicity restrictions.

An important goal for future work is to  apply the procedures to real-world data. 
However,  it is unreasonable to expect any non-simulated data to follow a recursive \ML\ model exactly, and the model should then be modified by adding appropriate noise terms. 
In particular we should not expect that we observe a minimal observed ratio more than twice, as we exploit in Algorithm~\ref{ch4:alg4}. It seems to be more reasonable to expect values close to each other. We therefore want to develop methods based on accumulation points.  It is hard to imagine noise models that would lead to simple exact likelihood analysis. One should then rather study the asymptotic precision of reasonable estimates and their behaviour under appropriate scaling, for example along the lines of \cite{davis:mccormick:89}. 

\section*{Acknowledgements}

We thank Justus Hartl for providing a first discussion about the different estimators suggested in this paper in his master's thesis. NG acknowledges  support by Deutsche Forschungsgemeinschaft (DFG) through  the TUM International Graduate School  of Science and Engineering (IGSSE). All authors benefited from financial support from the Alexander von Humboldt Stiftung.  

\begin{small}

\end{small}

\appendix
\section{Appendix: some technical proofs}
\subsection*{Proof of Lemma~\ref{ch4:lem}} 
\bproof 
First, define for $i\in V$ 
\begin{align*}
\Omega_{1/2}^{1,i}&:=\big\{ X_i= \bigvee_{k\in\pa(i)}b_{ki}X_k=\bigvee_{k\in\pa(i)}b^*_{ki}X_k\big\} , \quad \Omega_{1/2}^{2,i}:= \big\{X_i> \bigvee_{k\in\pa(i)}(b_{ki}\vee b^*_{ki})X_k\big\},\\
\Omega_i&:= \big\{ \bigvee_{j\in\An(i): b_{j i}= b^*_{j i} } b_{j i} Z_j> \bigvee_{j\in\an(i):  b_{j i}\neq b^*_{j i}}(  b_{j i} \vee  b^                                                                                                                                                                                                                                                             *_{j i})  Z_j\big\}. 
\end{align*}
The proof is by induction on the number of nodes of $\D$. For $d=1$ the statement is clear. Assume now that $\D=(V,E)$ has $d+1$ nodes and that the assertion  holds with respect to  \DAG s with at most $d$ nodes. Furthermore, assume without loss of generality that $d+1$ is a terminal node (i.e., $\des(d+1)=\emptyset$). Since $(X_1,\ldots, X_{d})$ 
follows a \mSEM\ on the \DAG\ $(\{1,\ldots,d\}, E\cap (\{1,\ldots,d\}\times \{1,\ldots,d\}))$ with \CM\ $ B=( b_{ij})_{d\times d}$ and  $ B^*=( b^*_{ij})_{d\times d}$ is the \ML\ coefficient matrix of a recursive \ML\ model on this \DAG\ as well, the induction hypothesis yields that
\begin{align}\label{ch4:lem:eq2}
\P(F\cap \{  \bfx \in A_{1/2}(B,B^*) \}) &= \P\big(F\cap   \bigcap_{i=1}^{d+1} \big(\Omega_{1/2}^{1,i} \cup \Omega_{1/2}^{2,i} \big)   \big)=\P\big(F\cap   \bigcap_{i=1}^{d} \Omega_i \cap \big( \Omega_{1/2}^{1,d+1} \cup \Omega_{1/2}^{2,d+1}\big) \big).
\end{align}
For every $i\in V$ we have by \eqref{ch4:ml-noise}  on $ \Omega_i$ that
\begin{align}\label{ch4:lem:eq1}
X_i=\bigvee_{j\in\An(i)} b_{j i}Z_j=\bigvee_{j\in\An(i)} b^*_{j i}Z_j.
\end{align} 
Noting from the proof of Theorem~4.2 of \cite{GK1} that 
\begin{align*} 
\bigvee_{k\in\pa(d+1)}b_{k,d+1}X_k= \bigvee_{k\in\pa(d+1)} b_{k,d+1} \bigvee_{j\in\An(k)}b_{j k}Z_j=\bigvee_{j\in\an(d+1)}b_{j,d+1} Z_j,
\end{align*}
we obtain from \eqref{ch4:lem:eq1} on $\bigcap_{i=1}^d  \Omega_i$, 
\begin{align*} 
\bigvee_{k\in\pa(d+1)}b^*_{k,d+1}X_k= \bigvee_{k\in\pa(d+1)} b^*_{k,d+1} \bigvee_{j\in\An(k)}b^*_{j k}Z_j =\bigvee_{j\in\an(i)}b^*_{j,d+1} Z_j.
\end{align*}
Thus, again by \eqref{ch4:ml-noise},
\begin{align*}
\bigcap_{i=1}^d \Omega_i\cap \Omega_{1/2}^{1,d+1}
&= \bigcap_{i=1}^d \Omega_i \cap \big\{ \bigvee_{j\in\An(d+1)} b_{j,d+1}Z_j=  \bigvee_{j\in\an(d+1)} b_{j,d+1} Z_j=\bigvee_{j\in\an(d+1)} b^*_{j,d+1} Z_j \big\},\\
\bigcap_{i=1}^d \Omega_i \cap \Omega_{1/2}^{2,d+1}&= \bigcap_{i=1}^d \Omega_i \cap \big\{ \bigvee_{j\in\An(d+1)} b_{j,d+1}Z_j>  \bigvee_{j\in\an(d+1)} ( b_{j,d+1}\vee b^*_{j,d+1} ) Z_j \big\}\\
&=\bigcap_{i=1}^d \Omega_i \cap \big\{  b_{j,d+1}Z_j>  \bigvee_{j\in\an(d+1)} ( b_{j,d+1}\vee b^*_{j,d+1} ) Z_j \big\}. 
\end{align*}
From \eqref{ch4:noiseNull} we then finally observe that  $\bigcap_{i=1}^{d} \Omega_i \cap \big( \Omega_{1/2}^{1,d+1} \cup \Omega_{1/2}^{2,d+1}\big)$ and $ \bigcap_{i=1}^{d} \Omega_i \cap \Omega_{d+1}$  
 only differ by a set of probability zero, and, hence,  \eqref{ch4:lem:eq3} follows from \eqref{ch4:lem:eq2}.
\eproof
\subsection*{Proof of Theorem~\ref{ch4:dens:general1}}
\bproof We must verify properties (A)--(C) of (\ref{eq:partition}).

\noindent(A)   Since $V$ is finite, it suffices to show for every $i\in V$,
\begin{align}
&P_{ B}\big(\big\{\bfsx\in\R_+^d: x_i< \bigvee_{k\in \pa(i)}  b_{ki}x_k\big\}\big)= \P\big(   X_{i}< \bigvee_{k\in \pa(i)}  b_{ki}   X_k\big)=0,\label{ch4:gmle:eq5}\\ \quad 
&P_{ B}\big(\{\bfsx\in\R_+^d:x_{i}  = \bigvee_{k\in \pa(i)} b^*_{ki}x_k> \bigvee_{k\in \pa(i)}  b_{k i}x_k \big \}\big)=\P\big(   X_{i}  = \bigvee_{k\in \pa(i)} b^*_{ki} X_k  > \bigvee_{k\in \pa(i)}  b_{ki}  X_k \big)= 0. \nonumber   
\end{align}
The former is immediate by \eqref{ch4:proof:dens}. By the same argument  we have for the latter,
\begin{align*}
0&\le \P\big( \bigvee_{k\in \pa(i)}  b_{ki} X_k \vee   Z_{i}= \bigvee_{k\in \pa(i)} b^*_{ki} X_k >\bigvee_{k\in \pa(i)}  b_{k i} X_k \big)
= \P\big(  Z_{i}= \bigvee_{k\in \pa(i)} b^*_{ki}  X_k >\bigvee_{k\in \pa(i)} b_{ki}  X_k \big)\\
&\le  \P\big(  Z_{i} =\bigvee_{k\in \pa(i)} b^*_{ki}\bigvee_{j\in\An(k)}  b_{j k}  Z_j \big)=0, 
\end{align*}
where we have used \eqref{ch4:ml-noise} and \eqref{ch4:noiseNull} for the last inequality and  equality, respectively. 
Thus we have verified (A). 

\noindent(B) Recall  that $P_B$ and $P_{B^*}$ share the same innovation vector when represented by a recursive ML model.   
Furthermore, note that the set $\Omega(B,B^*)$ from Lemma~\ref{ch4:lem} is a subset of $\bigcap_{i\in V}\big\{X_i=   \bigvee_{j\in\An(i): b_{j i}=b^*_{j i}}  b_{j i}Z_j\big\}$. 
We have $\Omega(B,B^*)=\Omega(B^*,B)$ and hence 
we obtain from \eqref{ch4:lem:eq3} for $A\in\borel(\R_+^d)$,
\begin{align*}
P_{B} (A\cap A_{1/2}(B,B^*)) &= \P(  \{ \bfx\in A\}\cap \Omega(B,B^*) )=\P\big( \big \{ \big(  \bigvee_{j\in\An(i): b_{j i}=b^*_{j i}}  b_{j i}Z_j,i\in V \big) \in A\big\}\cap \Omega(B,B^*) \big)\\
&=\P\big( \big \{ \big(  \bigvee_{j\in\An(i): b_{j i}=b^*_{j i}}  b^*_{j i}Z_j,i\in V \big) \in A\big\}\cap \Omega(B^*,B) \big)
= P_{B^*}(A\cap A_{1/2} (B,B^*)). 
\end{align*}

\noindent(C) We observe from the definition of  $A_0(B,B^*)$ and $A_{1/2}(B,B^*)$ that
\begin{align*}
\lefteqn{A_1( B,B^*)=\R_+^d\setminus \big(A_0(B,B^*) \cup  A_{1/2}(B,B^*)\big)} \\
&\subseteq\bigcup_{i\in V}\big[ \big\{ \bfsx\in\R_+^d:   \bigvee_{k\in \pa(i)} b^*_{ki}x_k > x_{i} \ge \bigvee_{k\in \pa(i)} b_{ki}x_k \big\}\cup \big\{ \bfsx\in\R_+^d: x_{i} = \bigvee_{k\in \pa(i)}  b_{ki}x_k >\bigvee_{k\in \pa(i)} b^*_{ki}x_k \big\}\big]\\
&\subseteq A_0(B^*,B). 
\end{align*}
Since  $A_0( B^*,B)$ is  a $P_{B^*}$-null set by (A),  this holds for the subset  $A_1( B,B^*)$ as well.
\eproof

\subsection*{Proof of Proposition~\ref{ch4:prop:rhoi}}
\begin{figure}
	\begin{center}	
	\begin{tabular}{ccc}
			\begin{tikzpicture}[scale=0.75]
			\coordinate (y) at (0,4);
			\coordinate (x) at (-0.2,0);
			\draw[<-] (y) node[above] {} -- (0,-0.2);
			\draw[->] (x)-- (4.5,0) node[below]{\begin{footnotesize}$x_i$\end{footnotesize}};
			\draw (-.1,0) -- (.1,0) node[left=4pt] {\begin{small}$\scriptstyle 0$\end{small}};
			\draw (-.1,3.5) -- (.1,3.5) node[left=4pt] {\begin{small}$\scriptstyle 1$\end{small}};
			\draw  (3,2.4)   node[right] {\begin{scriptsize}\color{OliveGreen}$F_{Z_i}$\end{scriptsize}};
			\draw   (0,4)   node[right] {\begin{tiny}\color{OliveGreen}$P^{i\mid\pa(i)}_B\big((0,x_i]\mid \bfsx_{\pa(i)}\big)$\end{tiny}};
			\draw[OliveGreen,line width=1pt,samples=100,domain=1.5:4.5] plot(\x,{3-3*exp(-0.75*\x)});
			\draw [OliveGreen,fill=OliveGreen]  (1.5,2.026) circle (2pt) node {};
			\draw [line width=0.5pt,,color=OliveGreen,dashed]  (1.5,1.95) -- (1.5,-5.1);

			\coordinate (y) at (0,-1);
			\coordinate (x) at (-0.2,-5);
			\draw[<-] (y) node[above] {} -- (0,-5.2);
			\draw[->] (x)-- (4.5,-5) node[below]{\begin{footnotesize}$x_i$\end{footnotesize}};
			\draw  (3,-2.6)   node[blue,right] 
			{\begin{scriptsize}$F_{Z_i}$\end{scriptsize}};
			\draw (-.1,-5) -- (.1,-5) node[left=4pt] {\begin{small}$\scriptstyle 0$\end{small}};
			\draw (-.1,-1.5) -- (.1,-1.5) node[left=4pt] {\begin{small}$\scriptstyle 1$\end{small}};
			\draw[blue,line width=1pt,samples=100,domain=2.5:4.5] plot(\x,{3-3*exp(-0.75*\x)-5});
			\draw [blue,fill=blue]  (2.5,-2.46) circle (2pt) node {};
			\draw [line width=0.5pt,blue,dashed]  (2.5, 2.6) -- (2.5,-5.3) ;
			\draw   (0,-1)   node[right] {\begin{tiny}\color{blue}$P^{i\mid\pa(i)}_{B^*}\big((0,x_i]\mid \bfsx_{\pa(i)}\big)$\end{tiny}};
			\draw   (1,-5.13)   node[below] {\begin{tiny}
				\color{OliveGreen}  $\bigvee\limits_{k\in\pa(i)} b_{ki}x_k$
				\end{tiny}};
			\draw   (3,-5)   node[blue,below=2pt] {\begin{tiny}
				$\bigvee\limits_{k\in\pa(i)} b^*_{ki}x_k$
				\end{tiny}};
			\end{tikzpicture}
			&
			
			\begin{tikzpicture}[scale=0.75]
			\coordinate (y) at (0,4);
			\coordinate (x) at (-0.2,0);
			\draw[<-] (y) node[above] {} -- (0,-0.2);
			\draw[->] (x)-- (4.5,0) node[below]{\begin{footnotesize}
				$x_i$\end{footnotesize}};
			\draw  (3,2.4)   node[right,OliveGreen] {\begin{scriptsize}$F_{Z_i}$\end{scriptsize}};
			\draw (-.1,0) -- (.1,0) node[left=4pt] {\begin{small}$\scriptstyle 0$\end{small}};
			\draw (-.1,3.5) -- (.1,3.5) node[left=4pt] {\begin{small}$\scriptstyle 1$\end{small}};
			\draw[OliveGreen,line width=1pt,samples=100,domain=2:4.5] plot(\x,{3-3*exp(-0.75*\x)});
			\draw [OliveGreen,fill=OliveGreen]  (2,2.33) circle (2pt) node {};
			\draw [line width=0.5pt,,color=OliveGreen,dashed]  (2,2.25) --  (2,-5.1) ; 
			\draw   (0,4)   node[right] {\begin{tiny}\color{OliveGreen}$P^{i\mid\pa(i)}_B\big((0,x_i]\mid \bfsx_{\pa(i)}\big)$\end{tiny}};
			
			\coordinate (y) at (0,-1);
			\coordinate (x) at (-0.2,-5);
			\draw[<-] (y) node[above] {} -- (0,-5.2);
			\draw[->] (x)-- (4.5,-5) node[below]{\begin{footnotesize}$x_i$\end{footnotesize}};
			\draw (-.1,-5) -- (.1,-5) node[left=4pt] {\begin{small}$\scriptstyle 0$\end{small}};
			\draw  (3,-2.6)   node[blue,right] 
			{\begin{scriptsize}$F_{Z_i}$\end{scriptsize}};
			\draw (-.1,-1.5) -- (.1,-1.5) node[left=4pt] {\begin{small}$\scriptstyle 1$\end{small}};
			
			\draw[blue,line width=1pt,samples=100,domain=2:4.5] plot(\x,{3-3*exp(-0.75*\x)-5});
			\draw [blue,fill=blue]  (2,-2.67) circle (2pt) node {};
			\draw   (0,-1.0)   node[right] {\begin{tiny}\color{blue}$P^{i\mid\pa(i)}_{B^*}\big((0,x_i]\mid \bfsx_{\pa(i)}\big)$\end{tiny}};
			\draw   (2.15,-5)   node[below=2pt] {\begin{tiny}
				$\color{OliveGreen}\bigvee\limits_{k\in\pa(i)} b_{ki}x_k={\color{blue}\bigvee\limits_{k\in\pa(i)} b^*_{ki}x_k}$\end{tiny}};
			\end{tikzpicture}
			&
			\begin{tikzpicture}[scale=0.75]
			\coordinate (y) at (0,4);
			\coordinate (x) at (-0.2,0);
			\draw[<-] (y) node[above] {} -- (0,-0.2);
			\draw[->] (x)-- (4.5,0) node[below]{\begin{footnotesize}$x_i$\end{footnotesize}};
			\draw (-.1,0) -- (.1,0) node[left=4pt] {\begin{small}$\scriptstyle 0$\end{small}};
			\draw  (3,2.4)   node[right] {\begin{scriptsize}\color{OliveGreen}$F_{Z_i}$\end{scriptsize}};
			\draw (-.1,3.5) -- (.1,3.5) node[left=4pt] {\begin{small}$\scriptstyle 1$\end{small}};
			\draw[OliveGreen,line width=1pt,samples=100,domain=2.5:4.5] plot(\x,{3-3*exp(-0.75*\x)});
			\draw [OliveGreen,fill=OliveGreen]  (2.5,2.540) circle (2pt) node {};
			\draw [line width=0.5pt,dashed,OliveGreen]  (2.5,2.25) -- (2.5,-5.2); 
			\draw   (0,4)   node[right] {\begin{tiny}\color{OliveGreen}$P^{i\mid\pa(i)}_B\big((0,x_i]\mid \bfsx_{\pa(i)}\big)$\end{tiny}};
			
			\coordinate (y) at (0,-1);
			\coordinate (x) at (-0.2,-5);
			\draw[<-] (y) node[above] {} -- (0,-5.2);
			\draw[->] (x)-- (4.5,-5) node[below]{\begin{footnotesize}$x_i$\end{footnotesize}};
			\draw (-.1,-5) -- (.1,-5) node[left=4pt] {\begin{small}$\scriptstyle 0$\end{small}};
			\draw  (3,-2.6)   node[blue,right]
			{\begin{scriptsize}$F_{Z_i}$\end{scriptsize}}; 
			\draw (-.1,-2) -- (.1,-2) node[left=4pt] {\begin{small}$\scriptstyle 1$\end{small}};
			\draw[blue,line width=1pt,samples=100,domain=1.5:4.5] plot(\x,{3-3*exp(-0.75*\x)-5});
			\draw [blue,fill=blue]  (1.5,-2.974) circle (2pt) node {};
			\draw [line width=0.5pt,blue,dashed] (1.5,-2.9) -- (1.5, -5.1) ;
			\draw   (0,-1)   node[right] {\begin{tiny}\color{blue}$P^{i\mid\pa(i)}_{B^*}\big((0,x_i]\mid \bfsx_{\pa(i)}\big)$\end{tiny}};
			\draw  (1,-5)   node[blue,below] {\begin{tiny}
				$\bigvee\limits_{k\in\pa(i)} b^*_{ki}x_k$
				\end{tiny}};
			\draw   (3,-5)   node[below=2pt] {\begin{tiny}\color{OliveGreen}$\bigvee\limits_{k\in\pa(i)} b_{ki}x_k$ \end{tiny}};
			\end{tikzpicture}
		\end{tabular}
		\caption{The conditional distribution functions from the proof of Proposition~\ref{ch4:prop:rhoi}(c).\label{ch4:fig:lemA4}
		}
	\end{center}
\end{figure}
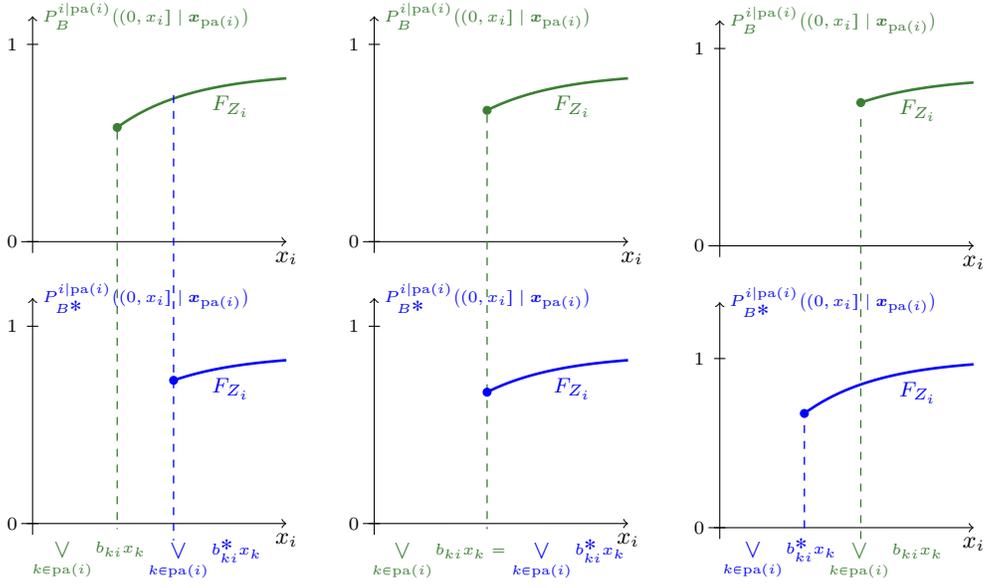

\bproof
Denoting by $A^i_0(B_i,B^*_i)$, $A^i_{1/2}(B_i,B^*_i)$, $A^i_1(B_i,B^*_i)$  the sets defining $\rho_i(\cdot, B_i,B^*_i)$, we have for the corresponding sets of $\rho$,
\begin{align*}
A_0(B,B^*)&= \bigcup_{i\in V} \big\{ \bfsx\in\R_+^d:  \bfsx_{\Pa(i)}\in A_0^i(B_i,B^*_i)\big\},\\
A_{1/2}(B,B^*)&= \bigcap_{i\in V} \big\{ \bfsx\in\R_+^d:  \bfsx_{\Pa(i)}\in A_{1/2}^i(B_i,B^*_i)\big\},\\
A_{1}(B,B^*)&= \bigcap_{i\in V} \big\{ \bfsx\in\R_+^d:  \bfsx_{\Pa(i)}\in A_{1/2}^i(B_i,B^*_i)\cup A_{1}^i(B_i,B^*_i)\big\} \cap  \big[ \R_+^d \setminus A_{1/2}(B_i,B^*_i)\big].
\end{align*}
From this we obtain (a) and (b).  Now, to see (c) we reason as follows:
\begin{align*} 
P_B^{i\mid \pa(i)}\big((0,x_i]\mid \bfsx_{\pa(i)}\big)&= F_{Z_i} (x_i) \bone_{[\bigvee_{k\in\pa(i)} b_{ki}x_k ,\infty)}(x_i), \quad \bfsx_{\Pa(i)}\in \R_+^{\vert \Pa(i)\vert},
\end{align*}
is a regular conditional distribution function of $X_i$ given $\bfx_{\pa(i)}$. To see this, use \eqref{ch4:proof:dens} and the independence of the innovations to obtain  
\begin{align*}
P_B^{i\mid \pa(i)}\big((0,x_i]\mid \bfsx_{\pa(i)}\big)&= \P(  X_i \le x_i\mid \bfx_{\pa(i)}=\bfsx_{\pa(i)} )\\
&= \P\big(   \bigvee_{k\in\pa(i)} b_{ki}X_k\vee Z_i \le x_i\mid \bfx_{\pa(i)}=\bfsx_{\pa(i)}\big )\\
&= F_{Z_i} (x_i) \bone_{[\bigvee_{k\in\pa(i)} b_{ki}x_k  ,\infty)}(x_i). 
\end{align*}
Since $\bfx$ and $\bfx^*$ share  the same innovation vector,  we have
\begin{align*} 
P_{B^*}^{i\mid \pa(i)}\big((0,x_i]\mid \bfsx_{\pa(i)}\big)&= F_{Z_i} (x_i) \bone_{[\bigvee_{k\in\pa(i)} b^*_{ki}x_k  ,\infty)}(x_i), \quad \bfsx_{\Pa(i)}\in \R_+^{\vert \Pa(i)\vert},
\end{align*} 
is a regular conditional distribution function of $X^*_i$ given $\bfx^*_{\pa(i)}$.  Figure~\ref{ch4:fig:lemA4} depicts  the two  conditional distribution functions for the three possible orders between $\bigvee_{k\in\pa(i)} b_{ki}x_k$  and $\bigvee_{k\in\pa(i)} b^*_{ki}x_k$. It then suffices to show for all $\bfsx_{\pa(i)}\in\R_+^{\vert\pa(i)\vert}$ and $y\in\R_+$, 
\begin{align*}
P^{i\mid\pa(i)}_{B}\big((0,y]\mid \bfsx_{\pa(i)}\big)=\int_{(0,y]} \rho_i(\bfsx_{\Pa(i)}, B_i,B^*_i) \big(P^{i\mid \pa(i)}_{B}+P^{i\mid \pa(i)}_{B^*}\big)(dx_i\mid \bfsx_{\pa(i)}),
\end{align*}
and for this again by definition of $\rho_i$ (cf. \eqref{ch4:dens:gen} and the related discussion) that
\begin{align*}
&P_B^{i\mid \pa(i)}\big((0,y]\cap \big(0, \bigvee_{k\in\pa(i)}b_{ki}x_k\big) \mid \bfsx_{\pa(i)} \big)=0,\\  
&P_B^{i\mid \pa(i)}\big((0,y]\cap \big\{\bigvee_{k\in\pa(i)}b^*_{ki}x_k \big\} \mid \bfsx_{\pa(i)} \big)=0 \quad \text{if $\bigvee_{k\in\pa(i)}b^*_{ki}x_k>\bigvee_{k\in\pa(i)}b_{ki}x_k$,}\\
&P_B^{i\mid \pa(i)}\big((0,y]\cap \big\{ \bigvee_{k\in\pa(i)} b_{ki}x_k\big\}\mid \bfsx_{\pa(i)}\big)=P_{B^*}^{i\mid \pa(i)}\big((0,y]\cap \big\{ \bigvee_{k\in\pa(i)} b_{ki}x_k\big\}\mid \bfsx_{\pa(i)}\big) \\
&\hspace*{19.05em} \text{if $\bigvee_{k\in\pa(i)}b^*_{ki}x_k=\bigvee_{k\in\pa(i)}b_{ki}x_k$,}\\ 
&P_B^{i\mid \pa(i)}\big((0,y]\cap \big(\bigvee_{k\in\pa(i)}(b_{ki}\vee b^*_{ki})x_k, \infty\big) \mid \bfsx_{\pa(i)}\big)=P_{B^*}^{i\mid \pa(i)}\big((0,y]\cap \big(\bigvee_{k\in\pa(i)}(b_{ki}\vee b^*_{ki})x_k, \infty\big) \mid \bfsx_{\pa(i)}\big). 
\end{align*}
Since $F_{Z_i}$ is atom-free, 
 this can be read directly from Figure~\ref{ch4:fig:lemA4}. 
\eproof
\end{document}